\documentclass[9pt]{article}
\usepackage{graphicx}
\usepackage{float}
\usepackage{amsmath}
\usepackage{epsfig}
\usepackage{bm}
\usepackage{url}
\usepackage[dvipsnames]{xcolor}
\usepackage{amsmath}
\usepackage{amsthm}
\usepackage{graphicx}
\usepackage{amssymb}
\usepackage{subfigure}
\usepackage[a4paper, total={6in, 8in}]{geometry}

\newtheorem*{remark}{Remark}
\newcommand{\indep}{\perp\!\!\!\!\perp} 
\DeclareUnicodeCharacter{02BC}{'}

\title{\textbf{Detection of Undeclared EV Charging Events in a Green Energy Certification Scheme  }}
\author{\small Luca Domenico Loiacono$^{1,2}$, Anthony Quinn$^{1,3}$, Emanuele Crisostomi$^{2}$, Robert Shorten$^{1}$}

\date{}

\begin{document}

\maketitle

\begin{center}
    \small $^{1}$Dyson School of Design Engineering, Imperial College London \\
    \small $^{2}$Department of Energy, Systems, Territory and Constructions Engineering, University of Pisa \\
    \small $^{3}$Department of Electronic and Electrical Engineering, Trinity College Dublin \\
\end{center} 

\begin{abstract}

\noindent
The green potential of electric vehicles (EVs) can be fully realized only if  their batteries are charged using energy generated from renewable (i.e.\ green) sources. For logistic or economic reasons, however, EV drivers may be tempted to avoid charging stations certified as providing green energy, instead opting for conventional ones, where  only a fraction of the available energy is green. This behaviour may slow down the achievement of decarbonisation targets of the road transport sector.  In this paper, we use GPS data to infer whether an undeclared charging event has occurred. Specifically, we construct a Bayesian hypothesis test for the charging behaviour of the EV. Extensive simulations are carried out for an area of London, using the mobility simulator, SUMO, and exploring various  operating conditions. Excellent detection rates for undeclared charging events are reported. We explain how the algorithm can serve as the basis for an incentivization scheme, encouraging  compliance by drivers with green charging policies.



\end{abstract}




\section{Introduction}
\label{sec:intro}


Electric vehicles (EVs) are emerging as pivotal solutions to mitigate urban air pollution, particularly when they replace conventional internal combustion engine (ICE) vehicles. According to the IEA's  2023 Global Electric Vehicle Outlook~\cite{iea2023global}, the global adoption of EVs is accelerating, and, by 2030, under current policies, EVs are expected to comprise about one-third of cars in China and approximately one-fifth in both the United States and the European Union. 

The benefits of EVs---in terms of reduced tail-pipe emissions---are undeniable,  particularly for people living in urban centres. 
However, the full potential of EVs can be realised only if renewable energy sources are used in EV charging. Only then will the uptake of EVs  actually contribute to the  significant decarbonisation of the future road transport sector, while providing co-benefits in terms of reduced air pollution\footnote{\url{https://www.eea.europa.eu/publications/electric-vehicles-and-the-energy} \\
This work has been submitted to the IEEE for possible publication. Copyright may be transferred without notice, after which this version may no longer be accessible.}. 

A   failure to guarantee the utilization of renewable energy in EV charging leads to the \textit{off-shoring} of emissions \cite{urry2014offshoring}.
Off-shoring refers to the shifting of pollution from an EV's location in a {\em source} country or region to a different---often geographically far removed--- {\em host\/} country or region where the electric energy used to charge its battery is actually generated. The negative impact of this practice becomes particularly significant if the electric energy is derived from fossil fuels (e.g.\ coal). Typically, such (polluting) energy production is off-shored to hosts  with laxer environmental regulations. Therefore, evaluating the environmental and social justice impacts of EVs is a complex issue that involves understanding how environmental policies and industrial practices affect pollution distribution. Unilateral environmental policies in source countries  can reduce {\em local\/} emissions by moving emission-intensive production abroad. However, this can lead to an increase in {\em global\/} emissions if the host countries have less stringent regulations, such as a heavy reliance on fossil fuels for electricity generation~\cite{jakob2021carbon}. 

Central to these concerns is whether EVs actually  provide any carbon footprint advantage over ICE vehicles, when the EV's production and end-of-life phases are accounted for. Analysis\footnote{\url{www.polestar.com/dato-assets/11286/1630409045polestarlcarapportprintkorr11210831.pdf}} of Polestar 2 EV variants, for instance, has shown that---in order to break even against the Volvo XC40 ICE vehicle---their required use phase (i.e.\ distance travelled) varies enormously depending of the combination of EV variant and the electricity mix they use. This ranges from 40,000 km for the Polestar 2 standard range single motor EV using wind-generated energy, to 110,000 km for the long range dual motor variant using a global electricity mix.    

To address these concerns about EV energy production---offshoring, off-setting, and the electricity mix---and to counter {\em green hypocrisy} among EV stakeholders, who appear to be environmentally friendly without actually being so, it is crucial to ensure that the energy used for EV charging truly has a green origin, i.e.\ that it comes from renewable energy sources (RES)~\cite{will2024can}. Various methodologies are used to guarantee that charging stations only use energy directly generated from RES to charge EVs, or, more often, that the energy used to charge EVs is offset from generation from RES  \cite{barman2023renewable,chang2021coordinate,mozafar2017simultaneous,tirunagari2022reaping}. We refer to these methodologies collectively as the {\em certification of green energy}. In this paper, we envisage that {\em certified charging stations\/}  adhere to a particular  green-energy certification scheme, and can, accordingly, incentivize EV drivers/owners (forthwith referred to simply as the driver) who decide to charge at these stations. 



For the EV driver,  charging at certified charging stations  may be  expensive or inconvenient when compared to uncertified charging with a generic energy mix, such as in the case of  home charging~\cite{allcost2024}. Accordingly, they may be induced to {\em hide\/} such a charging event, in order to avoid a notional penalty associated with uncertified charging, or, alternatively, in order to claim a bonus from a notional  incentivization scheme designed to encourage certified charging events. This paper is concerned with detecting whether such  {\em undeclared\/} (i.e.\ hidden) charging events have actually occurred between consecutive declared and certified charging events\footnote{While the most probable reason for a driver to hide a charging event is because it occurred at an uncertified charging point, our algorithm detects only whether an undeclared charging event occurred, and not whether the charge was certified. The idea here  is that  a {\em necessary\/} condition of compliance with an incentivization scheme for certified charging is that every charging event be declared.}.  

Our design of a detection algorithm for undeclared charging events  is based on a key assumption to ensure reliability; i.e.\  that all data available from the EV, including direct measurements of the state-of-charge (SoC), cannot be used to evaluate whether 
undeclared charging events have taken place, for the obvious reason that an  EV driver who wishes to hide a charging event may be induced also to interfere with the EV's data. Instead, we assume that a GPS device is  installed in the EV (as is usually the case  for insurance reasons), and we assume that the GPS data are reliable. During a certified charging event, the GPS device communicates to our algorithm the sequence of trips---specifically the time-stamped sequences of speeds and altitudes of the EV during those trips---since the last certified charging event.

The detection algorithm which we develop in this paper is Bayesian, in the sense that it elicits prior probability models for unobserved states, conditions on observations and known states, and constructs a predictive probability model. The latter yields the (detection) probability of an undeclared charging event (or events\footnote{We do not infer the number of such events, only that one or more has/have occurred.}) having occurred since the last certified charging event (which we will call the {\em certified interval} between charging). As already stated,  we process the  GPS data (notably the speed data),   and other known parameters of the EV (e.g.\ mass, battery size, type of heating/cooling device) to predict the  actual---but uncertain---energy consumption during the certified interval. Significant (i.e.\ improbable) deviations of the declared energy consumption from the predictions of the predictive probability model can be flagged, yielding a probability that an undeclared charging event occurred during the certified  interval.  We report extensive simulations using the mobility simulator, SUMO~\cite{behrisch2011sumo}. We show that the proposed methodology performs well in different seasons of the year (when the uncertainty of using auxiliary heating/cooling devices plays an important role in energy consumption), achieving true negative  (i.e.\ detecting that a compliant driver did not hide a charging event) rates  between $98.4\%$ and $99.9\%$. The true positive rates (i.e.\ detecting that non-compliant drivers did hide a charging event) reach $100\%$ in the case where the undeclared charge is  $50\%$ of the size of the battery. The details are provided in Section~\ref{sec:sims}.

As regulatory focus on off-shoring and energy sourcing intensifies, future policies are expected to consider not just the type of energy consumed by a vehicle, but also the source of that energy \cite{mckinsey2021}. Algorithms such as the one developed in this paper may serve as the enabling mathematical basis to ensure compliance with continuous green charging practices, and with evolving regulatory expectations.


\subsection{State of the art}
\label{sec:SOTA}
The main objective of this paper is to develop a Bayesian methodology which tests whether drivers 
are confining their charging events to certified green charging stations, and this can form the basis for a down-stream  incentivization scheme.   The idea is to test whether the data recorded by the GPS device during a certified  interval are consistent  with the change (i.e.\ differential) in the state-of-charge (SoC)---which is an energy quantity---
recorded by the consecutive certified charging stations. To the best of our knowledge, this is the first time that this problem has been formulated and discussed. Nevertheless, it is related to the problem of estimating the energy consumed by an EV over a specific route, for which many physics-based and data-driven methodologies have already been proposed in the literature, which we will now review. 

It is known that many factors influence EV energy consumption and can cause  significant variability over a specified GPS-tracked route, leading to  uncertainty in predicting the actual energy consumption. These factors include travel-related aspects, such as distance travelled, speed, acceleration, cumulative change in altitude, and cruising time \cite{brady2016development}; environmental factors such as ambient temperature \cite{jaguemont2016comprehensive,piao2022challenges}; visibility and wind effects \cite{donkers2020influence}; vehicle-related parameters including mass,  and auxiliary energy loads due to heating, ventilation, and air-conditioning (HVAC) \cite{kambly2014estimating,liu2018exploring,valentina2014hvac}; roadway topography such as slope \cite{liu2017impact};  traffic conditions \cite{madhusudhanan2020effect}; and driver-related aspects such as driving behaviours, car-following behaviour, and charging habits \cite{liu2016modelling,tang2015influences}. In addition, battery health and degradation effects due to ageing of the EV \cite{wang2021review,wu2022soc}---and even efficiency of regenerative braking \cite{braun2017influence,salari2023new}---are typical factors affecting EV energy efficiency. These considerations underscore the importance of comprehensively accounting for diverse variables to develop accurate models of EV energy consumption. While all the aforementioned factors do, indeed, influence the actual EV energy consumption, the most significant  are 
\begin{itemize}
\item[(i)] the net mass of the loaded EV, which is uncertain as it depends on the undeclared number and masses of the people; and 
\item[(ii)] the use of auxiliary power devices in the EV, most notably for heating and cooling purposes.
\end{itemize}
In our work, these two unobserved states provide the principal  sources of uncertainty about the actual energy consumed by a driver during the trips which they undertook in a certified  interval. This is  consistent with the energy consumption model (reviewed in Section~\ref{sec:physmodel}) developed in SUMO \cite{behrisch2011sumo}. SUMO is the popular road traffic simulator that we use in this paper (Section~\ref{sec:sims}) to model our road network and charging station infrastructure, as a means for validatation of our methodology \cite{kurczveil2014implementation}. 

In addition to the literature on estimation of energy consumption, there are related papers on the need to incentivize EV charging via renewable energy sources. A pioneering work in this direction was \cite{Sponge}, where the operation of EVs was orchestrated in order to maximise the consumption of energy generated from renewable sources. Recent advances in incentivizing EV charging with RES have led to innovative control algorithms designed to encourage accurate reporting of charging preferences (in terms of time and location) and optimize grid stability. For instance, 
\cite{kajanova2022social} introduced algorithms that incentivise truthful preference reporting among EV drivers, facilitating better alignment of charging patterns with renewable energy availability. Additionally,  \cite{soares2017dynamic} developed dynamic control algorithms that adjust incentives based on real-time data, in order to motivate EV drivers to adapt their charging schedules, thereby contributing to the management of grid stability and  to congestion avoidance. Similarly,  \cite{park2022electric} proposed financial rewards to encourage EV users to charge at specific times and locations, aiming to maximize the use of renewable energy and minimize shortages. 
These previous studies have principally  aimed to align EV charging with renewable energy production via financial incentives and real-time adjustments. In contrast, our approach is significantly different as it verifies adherence to incentive schemes by constructing a probability-based predictor of EV energy consumption during a certified  interval (i.e.\ between certified green charging events). This enables the detection of undeclared---and presumably uncertified---charging events, thereby encouraging constant compliance with certified green energy charging practices.

\subsection{Layout of the paper and notational conventions}
In Section~\ref{sec:physmodel}, we present our methodology for data-driven stochastic prediction of  the actual EV energy consumption during a certified interval.  This is used in Section~\ref{sec:hypo}  to compute a  hypothesis test (detector), which processes the recorded differential SoC during the certified interval into a probability that the driver has hidden (i.e.\ failed to declare) at least one charging event. Section~\ref{sec:sims} gives the details of the mobility simulator used for our validating case study, and summarizes the outcomes of the simulations. Finally, Section~\ref{sec:disc} concludes our manuscript by summarising our findings, noting the downstream role of the methodology in an incentivization scheme for green energy charging,  and outlining our current lines of research to extend the presented work.





\section{Energy consumption from an EV battery, $x_C$}
\label{sec:physmodel}
In what follows, we review a physical model (Section~\ref{Physical_Model}) which relates the energy consumed from  an EV battery  (i.e.\ drawn from the battery) during a trip (or trips) to the GPS-measurable properties of the trip(s) (Section~\ref{GPS}), and to   a set of known parameters (constants) associated with the specific EV. The model also depends on other unknown (i.e.\ unobserved) states of the EV, which may vary  during and between the trips. These are modelled stochastically in Section~\ref{sec:states}. This will enable construction of  a Bayesian predictive model of $x_C$ (Section~\ref{sec:hypo}), where $x_C$ is the total energy consumed from the battery during a certified  interval.

\subsection{Physical model}
\label{Physical_Model}

\noindent We adopt the physical model proposed in \cite{kurczveil2014implementation} and implemented in SUMO, in which the total energy consumed by  an EV during a trip 
is computed in a discrete-time manner. In what follows---and in the simulations in Section~\ref{sec:sims}---we assume a stepsize of 1 s, i.e.\ $t\in \{0,1,2,\ldots\}$, and all time-varying quantities (e.g.\ the speed, $v[t]$) are assumed to be constant during these 1 s time intervals. Where symbols are not defined in place below, they are {\em known\/} EV parameters---denoted by $\theta$---which are listed in Table~\ref{tab:vehicle_params}. 
\begin{enumerate}
    \item At (absolute) time $t\in\{0,1,\ldots\}$, the (instantaneous) energy of the EV is computed by summing its kinetic, potential, and rotational energy 
    components,
\begin{eqnarray}
            E_{veh}[t] & = & E_{kin}[t] \; + \; E_{pot}[t]  \;+ \; E_{rot,int}[t] \nonumber \\[4pt]
        &=&  \frac{1}{2} m_{total}[t] \times v^2[t] \; +\;  m_{total}[t] \times g \times h[t] \; +\;  \frac{1}{2} J_{int}\times v^2[t], \nonumber
        \label{First_Equation}
    \end{eqnarray}
where $m_{total}[t]$, $h[t]$ and $v[t]$ denote the total mass of the EV, its altitude and speed at time  $t$, respectively, and $g$ is the acceleration due to gravity.  We emphasize the discrete nature of the independent time variable, $t$, by adopting square brackets, $[\cdot]$.
    \item The energy required during the (next) time interval, $(t,t+1]$, is therefore
    \begin{equation}
      \Delta E_{req}[t]=E_{veh}[t+1]-E_{veh}[t].
    \end{equation}
   \item The energy consumed in time interval, $(t,t+1]$, can be calculated as follows:  
    \begin{equation}
        \Delta E_{cons}[t]=\Delta E_{req}[t]+\Delta E_{loss}[t].
    \label{eq:DelEcons}
    \end{equation}
    Here, the energy loss, $\Delta E_{loss}$, itself consists of three components,  
    \begin{eqnarray}
        \Delta E_{loss}[t]&=&\Delta E_{air}[t]+\Delta E_{roll}[t]+\Delta E_{aux}[t],\;\;\mathrm{where}
    \label{eq:Eloss}\\
        \Delta E_{air}[t]&=&\frac{1}{2} \times \rho_{air} \times A_{veh} \times c_w \times v^2[t] \times \Delta s[t], \nonumber\\
            \Delta E_{roll}[t]&=&c_{roll} \times m_{total}[t] \times g \times \Delta s[t],\label{eq:DelEroll}\\
            \Delta E_{aux}[t]&=&W[t],  \;\;\mathrm{and} \label{eq:DelEaux}\\
           \Delta s[t]&=&v[t]. \label{eq:DelEloss}
    \end{eqnarray}
    In (\ref{eq:Eloss}), $ \Delta E_{air} [t]$ and $\Delta E_{roll} [t]$ denote the energy loss in the EV due to air resistance and roll, respectively. They depend on   $\Delta s[t]$, i.e.\ the (absolute path) distance  travelled  in the 1 s interval, $(t,t+1]$. Also,   $\Delta E_{aux}[t]$ denotes the energy loss due to the use of auxiliary devices and services (e.g.\  cooling/heating), where $W[t]$ 
    is the  average power consumed by these  devices and services during the 1 s  interval.  Finally, $\rho_{air}$ denotes the density of air.
      \item 
      \(\Delta E_{cons}[t]\) may either be a positive or a negative quantity, since it represents the amount of energy the EV consumes due to its movement. It is negative, for instance,  when regenerative braking is used. The energy actually drawn from the battery in the 1 s interval, $(t,t+1]$, is denoted by $\Delta x[t]$. This depends on constant efficiency factors for the battery (Table~\ref{tab:vehicle_params}), i.e.\ $\eta_{prop} \in (0,1)$ if $\Delta E_{cons}[t]>0$ (the {\em propulsion\/} case), and  $\eta_{recup} \in (0,1)$ if $\Delta E_{cons}[t]<0$ (the {\em recuperation\/} case).  Specifically:
      \begin{eqnarray}
          \Delta x[t] &=& \Delta E_{cons}[t]\times \eta_{prop}, \;\; \mathrm{if} \;\; \Delta E_{cons}[t]>0, \nonumber\\
          \Delta x[t] &=& \Delta E_{cons}[t]\times \eta_{recup}, \;\; \mathrm{if} \;\; \Delta E_{cons}[t]<0. \label{Final_Equation}
      \end{eqnarray}
     Finally, we define $x_C[t]$ as  the {\em cumulative energy drawn from the battery\/} in a {\em net\/} interval of $t$ s. It is computed as the sum of the energies (\ref{Final_Equation}) drawn from the battery in each 1 s interval:
     \begin{equation}
     x_C[t]=x_C[t-1] +  \Delta x[t-1], \;\; t\in\{1,2,\ldots\},
     \label{eq:cumenergy}
     \end{equation}
initialized with $x_C[0]=0$.
     This unobserved (i.e.\ uncertain) time-varying quantity (i.e.\ state  process) is our principal quantity of interest, and  we will stochastically model and predict it in Section~\ref{sec:hypo}. 
     \end{enumerate} 

\begin{table}[H]
\centering
\begin{tabular}[t]{r|c|l}
\textbf{Parameter} & \textbf{Symbol} & \textbf{Units } \\ \hline
Battery Capacity & $E_{max}$  & kWh  \\ \hline
{Kerb Mass} & $m_{veh}$ & kg  \\ \hline
Front Surface Area & $A_{veh}$ & m$^2$  \\ \hline
Internal Moment Of Inertia & $J_{int}$ & 
$kg \cdot m^2$  \\ \hline
Radial Drag Coefficient & $c_{rad}$ & -  \\ \hline
Roll Drag Coefficient & $c_{roll}$ &   - \\ \hline
Air Drag Coefficient & $c_{w}$ &   - \\ \hline
Propulsion Efficiency & $\eta_{prop}$ & -   \\ \hline
Recuperation Efficiency & $\eta_{recup}$ & -   \\ 
\end{tabular}
\caption{EV parameters, $\theta$, used in the physical model (\ref{First_Equation})-(\ref{eq:cumenergy}).}
\label{tab:vehicle_params}
\end{table}

\subsection{GPS-measured  quantities, $v[t]$ and $h[t]$}
\label{GPS}
We have seen how $x_C[t]$---the cumulative energy drawn from the battery during a trip of duration $t$ s---depends on the sequentially measured (i.e.\ observed) speed and altitude of the EV, $v[t]$ and $h[t]$, respectively.  We have assumed that a GPS device  provides these measurements. An EV \textit{trip} refers to a distinct segment of driving activity, defined from the moment an EV 
is switched on
until it comes to a complete stop and is switched off.  Each trip is assigned a unique numeric label (Table~\ref{tab:GPS_dataframe}).  Note that $\Delta t =1$ s, as assumed in (\ref{eq:DelEloss}). Also, the discrete-time index, $t\in\{1,2,...\}$, is reset to $1$ in  (\ref{eq:cumenergy}) when the trip number increments, and  $x_c[0]$ is re-initialized with the cumulative value attained at the end of the previous trip.  In the sequel, ${\mathbf d}_G$ will denote the multi-trip GPS data record for the certified interval.
\begin{table}[H]
    \centering
    \begin{tabular}[t]{c|c|c|c}
    \hline
    \textbf{Trip Number} &  \textbf{GPS-time}  & \textbf{Speed} (\textit{v}) & \textbf{Altitude} (\textit{h}) \\ \hline
    1 &  07/07/24 10:23:59 & 3.4 \textit{m/s} & 0 \textit{m} \\
    1 &  07/07/24 10:24:00 & 1.9 \textit{m/s} & 2 \textit{m} \\
    1 &  07/07/24 10:24:01 & 10.4 \textit{m/s} & 4 \textit{m} \\ 
    1 &  07/07/24 10:24:02 & . & . \\
    1 &  07/07/24 10:24:03 & . & . \\ \hline
    2 &  08/07/24 19:03:14 & 1.1 \textit{m/s} & 30 \textit{m} \\ 
    2 &  08/07/24 19:03:15 & 5.7 \textit{m/s} & 33 \textit{m} \\
    2 &  08/07/24 19:03:16 & . & . \\
    2 &  08/07/24 19:03:17 & . & . \\ \hline
    
    \end{tabular}
    \caption{Extracts from a typical multi-trip GPS data record, ${\mathbf d}_G$, with $\Delta t =1 $ s.  }
    \label{tab:GPS_dataframe}
\end{table}

\subsection{Unobserved variables}
\label{sec:states}
Evaluation of the cumulative energy drawn from the battery, $x_C [t]$ (\ref{eq:cumenergy}),    also depends on other  unknown and unmeasured quantities, as follows:
\begin{itemize}
\item[(i)] $m_{total}[t]$, being the {\em total\/}  mass of the loaded EV at each measurement time, $t$ (\ref{First_Equation}), (\ref{eq:DelEroll});   
\item[(ii)] $W[t]$, the average power consumed by auxiliary services (such as onboard infotainment systems and air-conditioning), at each $t$ (\ref{eq:DelEaux}).
\end{itemize}
In the next two subsections,  we develop (static) probability models for these uncertain  processes.

\subsubsection{Total mass of the loaded EV, $m_{total}[t]$}
\label{sec:Tmass}
The total vehicular mass, $m_{total}[t]$, is given by the sum of the kerb mass, $m_{veh}$---which is an invariant  vehicular parameter (Table~\ref{tab:vehicle_params})---plus the net mass of people (driver plus passengers), $m_{peop}[t]$. In principle, one should also consider the mass of luggage, but  we neglect this for simplicity, as it is  small compared to that of the people. Importantly, we only model the  average mass of the people carried by the EV at {\em any\/} time during a trip, eliminating the dependence of $m_{peop}[t]$ on $t$; i.e.\ $m_{peop}[t] \equiv m_{peop}$. 

We adopt a Gaussian mixture model (GMM) as  the probability model for  $m_{peop}$, where each component models the different possible numbers of people typically carried by the EV, denoted by $N_p$ (to include the driver). The implied  hierarchical probability model is as follows:

\begin{itemize}
    \item \textbf{Prior  modelling of the number of people carried by the EV} 
    
    The  prior model for the number of people, 
    $ p_k \equiv \Pr[N_p = k],\;  k \in \left\{1,\ldots,5\right\}$.
This probability mass function  (pmf) is inferred using statistical data sourced from \cite{millard2011we} and Eurostat statistics \cite{CiteEursotatMobility}. From these sources, the expected number of people per  trip is found to be approximately $1.6$. The pmf is then constructed with this expected value as a constraint,  and assigning strictly decreasing probabilities. With five degrees-of-freedom but only two constraints, the pmf is under-determined. It is chosen as follows, to reflect the available data sources (above): 
\begin{table}[ht]
\centering
    \begin{tabular}[t]{r|l}
    \textbf{No.\ people, $N_p=k$} & \textbf{pmf}, $p_k$ \\
    \hline
        1 &  0.61 \\
        2 &  0.23 \\
        3 &  0.11 \\
        4 &  0.04 \\
        5 &  0.01 \\
    \end{tabular}
\end{table}   
    \item \textbf{Conditional Gaussian modelling of the mass of people}

    Given \( N_p = k \), the  mass of people carried by the EV, \( m_{peop}\), follows a Gaussian distribution:
    \begin{equation}
    \label{eq:massnor}
     \mathsf{F}(m_{peop} \mid N_p = k) \propto \mathcal{N}(m_{peop} \mid \mu_k, \sigma_k^2)
    \end{equation}
    with \(  m_{peop} \in \mathbb{R}^+ \). We choose the conditional mean to be   \( \mu_k = k \times 74 \),  and  the conditional standard deviation to be \( \sigma_k = 12 \sqrt{k} \). 
    The  $k=1$ standard interval (one person) is based on  data from \cite{walpole2012weight}.

    \item \textbf{Marginal modelling of  net mass (of people) carried by  EV}
    
    The overall distribution of mass (of people) carried by the EV is therefore the  mixture of components (\ref{eq:massnor}), weighted by the $N_p$-probabilities above:       \begin{equation}
    \mathsf{F}(m_{peop}) = \sum_{k=1}^{5} p_k\,\mathsf{F}(m_{peop} \mid N_p = k).
    \label{eq:GMM}
    \end{equation}
    This GMM is depicted  in Figure \ref{MixturePassengersMass}.
  \end{itemize}
\begin{figure}[H]
    \centering
    \caption{Gaussian mixture model (GMM), $\mathsf{F}(m_{peop})$, of the net mass of people typically carried by an EV during a trip. }
    \label{MixturePassengersMass}
\end{figure}

\subsubsection{Average power consumption of auxiliary devices, $W[t]$}
\label{sec:auxpower}
The energy consumed by auxiliary devices,  $E_{aux}$ (\ref{eq:Eloss}), can account for a significant proportion of  the total energy consumed by the EV during a trip, $E_{cons}$ (\ref{eq:DelEcons}). $E_{aux}$  is strongly dependent on the time of year, since  weather conditions strongly influence the utilization of HVAC  (Heating, Ventilation, and Air Conditioning) and  the lights \cite{al2022effects}, \cite{liu2018exploring}. For this reason, we {\em conditionally\/} model the unknown and unobserved auxiliary power, $W$ (\ref{eq:DelEloss}),  via the observed Bernoulli random variable (r.v.), $Y\in\{s,w\}$. These states denote  the case where the trip takes place during the warmer (summer\footnote{In this paper, {\em summer\/} might best be understood as {\em not winter}, as far as the London context of the SUMO simulation is concerned (Section~\ref{SUMO}). In that city, we associate summer with a low usage of any  auxiliary services, and winter with high usage of heating and lights. In other climates, of course, the high- and low-usage periods of the year may be delimited differently. }, $s$) 
or cooler (winter, $w$) periods of the year, respectively, a condition easily verified by the GPS time-stamps (Table~\ref{tab:GPS_dataframe}). 

\subsubsection*{Conditional Gamma Modelling of the Auxiliary Power, $W[t]$}

\begin{figure}[H]
    \centering
    \subfigure[]{\includegraphics[width=0.45\linewidth]{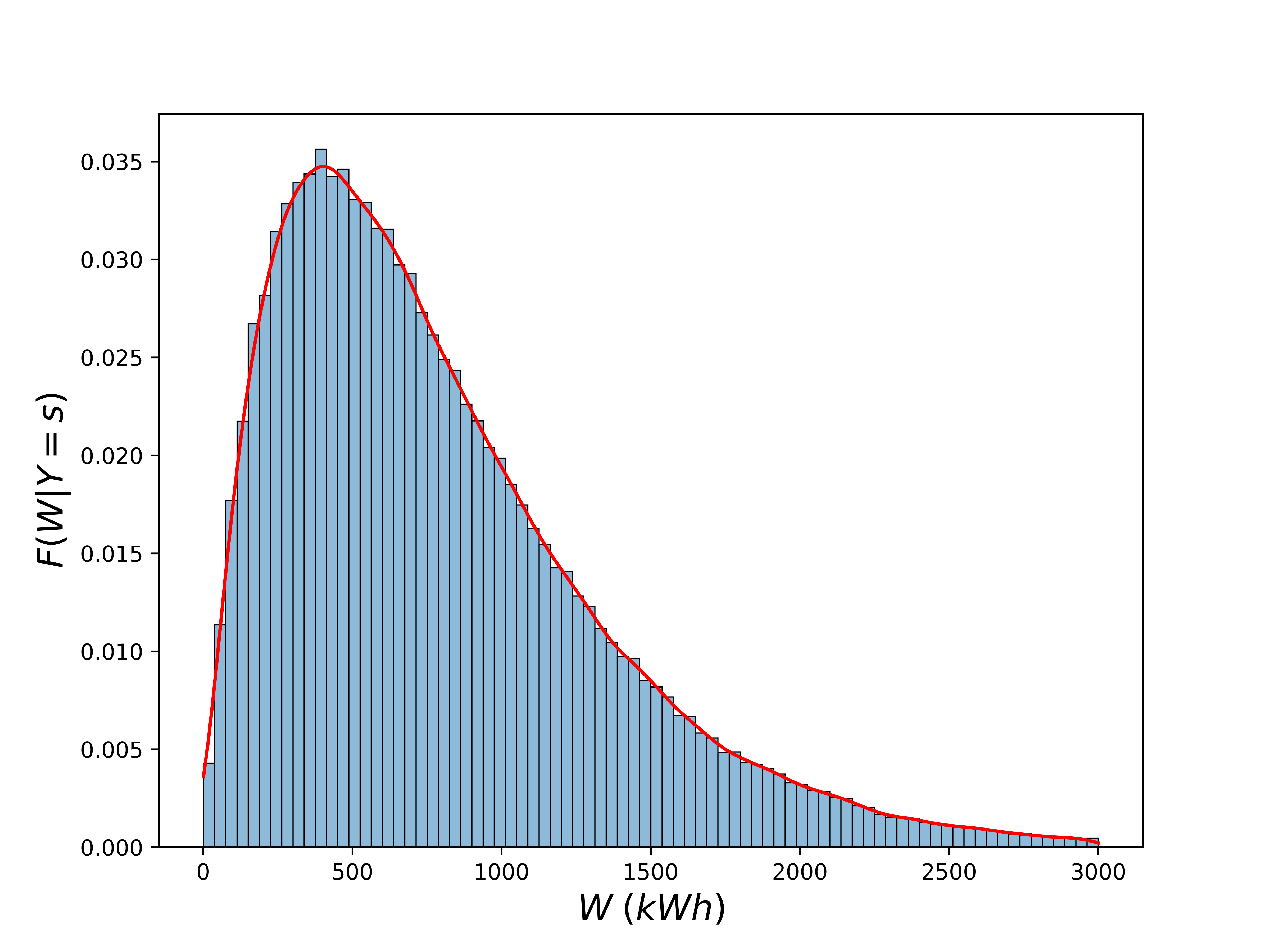}}
    \subfigure[]{\includegraphics[width=0.45\linewidth]{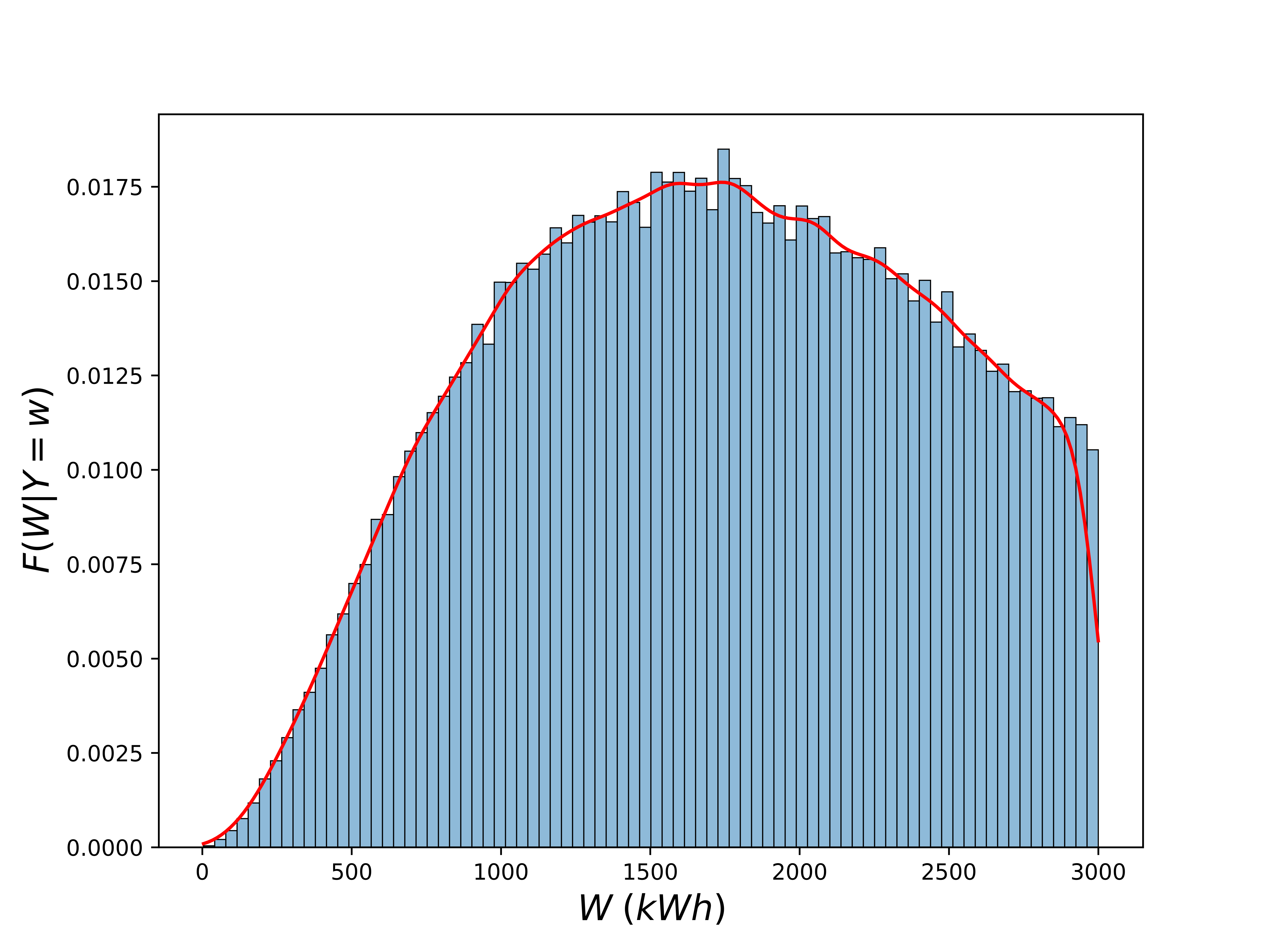}}
    \caption{Gamma models, $\mathsf{F}(W\mid Y=y)$, of the average power consumed by the auxiliary devices and services of an EV during a trip in (a) summer ($Y=s$), and (b) winter ($Y=w$).  }
    \label{fig:Gam}
\end{figure}

As in Section~\ref{sec:Tmass}, we facilitate the (conditional) stochastic modelling of the  auxiliary power process, $W[t]$, by modelling its {\em average\/} value, so that $W[t]\equiv W$ in (\ref{eq:DelEloss}). Since $W \in {\mathbb R}^+$, gamma modelling is appropriate~\cite{Bernardo:1319894}: 
\begin{equation}
\mathsf{F}(W|Y=y) \equiv
{\mathcal G}amma\left(k_y,\theta_y\right). 
\label{eq:Gam}
\end{equation}
Here, the $y\in\{s,w\}$-indexed parameters of the gamma model are its shape parameter, 
 $k_y > 0$, and its  scale parameter,
   $\theta_y > 0$, so that  
the conditional expected value (mean) and conditional variance of \( W \) are 
\begin{equation}
{\mathsf E}\left[W\middle| Y=y\right]={k}_y{\theta}_y,  \     \  \text{ Var}\left[W\middle| Y=y\right]={k}_y{\theta}_y^{2}.
\end{equation}
If these moments are available---for example via test data labelled by season---then the seasonally-conditioned gamma parameters, $k_y$ and $\theta_y$, can be estimated via moment-matching: 
\begin{equation}
    k_y \equiv \frac{{\mathsf E}[W| Y=y]^2}{\text{Var}[W| Y=y]}, \quad \theta_y \equiv \frac{\text{Var}[W| Y=y]}{{\mathsf E}[W| Y=y]}.
\label{eq:GamPar}
\end{equation}
Seasonal quantifiers of expected auxiliary power requirements, ${\mathsf E}\left[W\middle| Y=y\right]$, are available in the literature \cite{evtimov2017energy}, \cite{sun2021energy}, \cite{unni2021influence}, but we were not able to source the variances,  $\text{ Var}\left[W\middle| Y=y\right]$. Therefore, we adjusted the latter, constrained by the former (for each period of the year) to ensure a good fit of the published empirical data for $W$ (Figure~\ref{fig:Gam}), yielding the following conditional gamma parameters via  (\ref{eq:GamPar}):  
\begin{table}[ht]
\centering
\begin{tabular}[t]{lccc}
 & \textbf{Scale $\theta_y$} & \textbf{Shape ${k}_y$} & \( E[W| Y=y] \)  \\
\hline
        \textbf{Winter ($Y=w$)} & 800 & 3 & 2400 W\\
        \textbf{Summer ($Y=s$)} & 400 & 2 & 800 W\\
\end{tabular}
\end{table}%

\section{The probability, $ \Pr[H_1 \mid {\mathbf d}_G, \theta, y, x_D]$,   of undeclared charging during a certified  interval  }
\label{sec:hypo}

We now develop the main result of the paper, which is a hypothesis test for compliance with a notional green-energy charging scheme. Specifically, we  evaluate the probability  of the proposition, $H_1\equiv$ {\em the EV was charged during a specific (i.e.\ the last\footnote{We envisage that the test is implemented at every charging station belonging to the specific green energy scheme offering the incentives or applying the penalties of that scheme.}) certified  interval, without such charging event(s) having been declared}. This hypothesis test, yielding 
\begin{equation}
0<\Pr[H_1 \mid {\mathbf d}_G, \theta, y, x_D]<1,
\label{eq:hyptest}
\end{equation}
is to be implemented at the current certified charging station in order to detect whether undeclared charging occurred (for example, at home) at least once since the last certified charging event.  Since any undeclared charging event(s) cannot be certified by the green-energy scheme,  this hypothesis test (\ref{eq:hyptest}) can be used to withhold an incentive from---or even apply a penalty to---the driver under that scheme. 

When the EV connects at a certified charging station, we assume that the $H_1$-testing algorithm  can access the multi-trip GPS data, ${\mathbf d}_G$ (Table~\ref{tab:GPS_dataframe}), collected from all the trips that occurred during  the certified  interval. This, of course, means that  the season, $Y=y\in\{s,w\}$ (\ref{eq:Gam}), is also known. Furthermore, we assume that 
the specific EV model is recognized, and so its parameters, $\theta$ (Table~\ref{tab:vehicle_params}), can be looked up by the algorithm. 

Finally, and crucially, the charging station senses the  (absolute) SoC (a measure of energy, in kWh) of the EV immediately  after it is plugged in. We denote this by $x_1$. The algorithm must also have access to the SoC of the EV which was recorded immediately after charging at the previous certified charging station, i.e.\ $x_0$. The certified {\em differential SoC\/}, $x_D \equiv x_0 - x_1$, is therefore  observed (i.e.\ it appears in the condition of the {\em a posteriori\/} probability (\ref{eq:hyptest})). The purpose of the $H_1$-testing algorithm is to compare $x_D$ to the uncertain cumulative energy actually drawn from the battery during that interval, $x_C [T_c] \equiv x_C$ (\ref{eq:cumenergy}) (where $T_c$ denotes the cumulative time of all the trips taken during the certified interval\footnote{We will refer to $T_c$ simply as the {\em certified interval duration}, with the understanding that the time counter is reset at the beginning of each trip, as explained at the beginning of Section~\ref{GPS}.}). This probabilistic comparison  is inferentially equivalent to the soft (i.e.\ Bayesian) classification of $x_D$. Therefore, we must  first construct the distribution of  uncertain $x_C$ by processing (i.e.\ conditioning on) the entire knowledge base at the start of charging, this being ${\mathbf d}_G$, $\theta$, $Y=y$, along with  the uncertainty models (i.e.\ probability distributions), $\mathsf{F}(m_{peop})$ (\ref{eq:GMM}) and $\mathsf{F}(W|Y=y)$, $y\in\{s,w\}$ (\ref{eq:Gam}). We now present the various steps involved in constructing this test.     

\subsection{The battery energy conservation equation}
\label{sec:ChargeCons}

The observed SoC of the EV battery immediately after the {\em previously\/} certified charging event, $x_0$, and immediately before charging commences during the {\em current\/} certified charging event, $x_1$, are related to the (unobserved and uncertain) cumulative energy drawn from the battery during that interval, $x_c$, via the following energy conservation equation\footnote{\label{foot:trivial} Detection of undeclared charging is trivial when $x_D \leq 0$, since, then, $x_U \geq x_C > 0$.}:
\begin{equation}
        x_0-x_1 \equiv x_D \equiv x_{C} - x_{U} >0. 
\label{eq:conser}
\end{equation}
Here, \(x_U\) denotes the total amount of energy provided to the battery by undeclared charging events during the certified interval. We make the following key simplifying stochastic modelling assumption in respect of the r.v.s, $x_C$ and $x_U$: 


\[
x_{C} \indep x_{U},
\]
i.e.\ the amount of undeclared charging is independent of the cumulative energy drawn from the battery. We will use this assumption to elicit the distributions of $x_C$ and $x_U$ {\em independently}, in the next two subsections, respectively.

\subsection{Empirical probability model of $x_C$, the cumulative energy drawn from the battery}
\begin{figure}[H]
    \centering
    \subfigure[]{\includegraphics[width=0.45\linewidth]{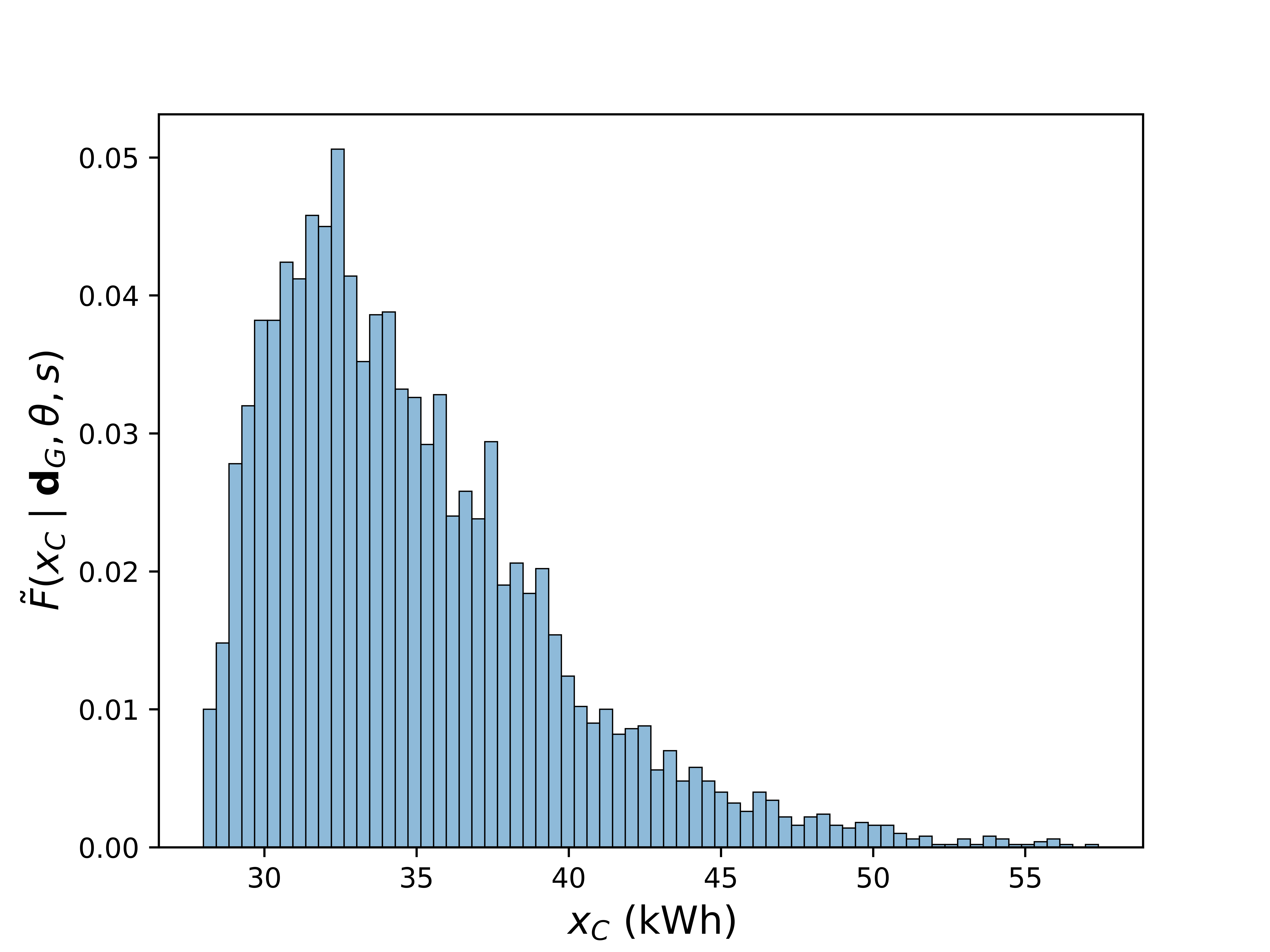}}
    \subfigure[]{\includegraphics[width=0.45\linewidth]{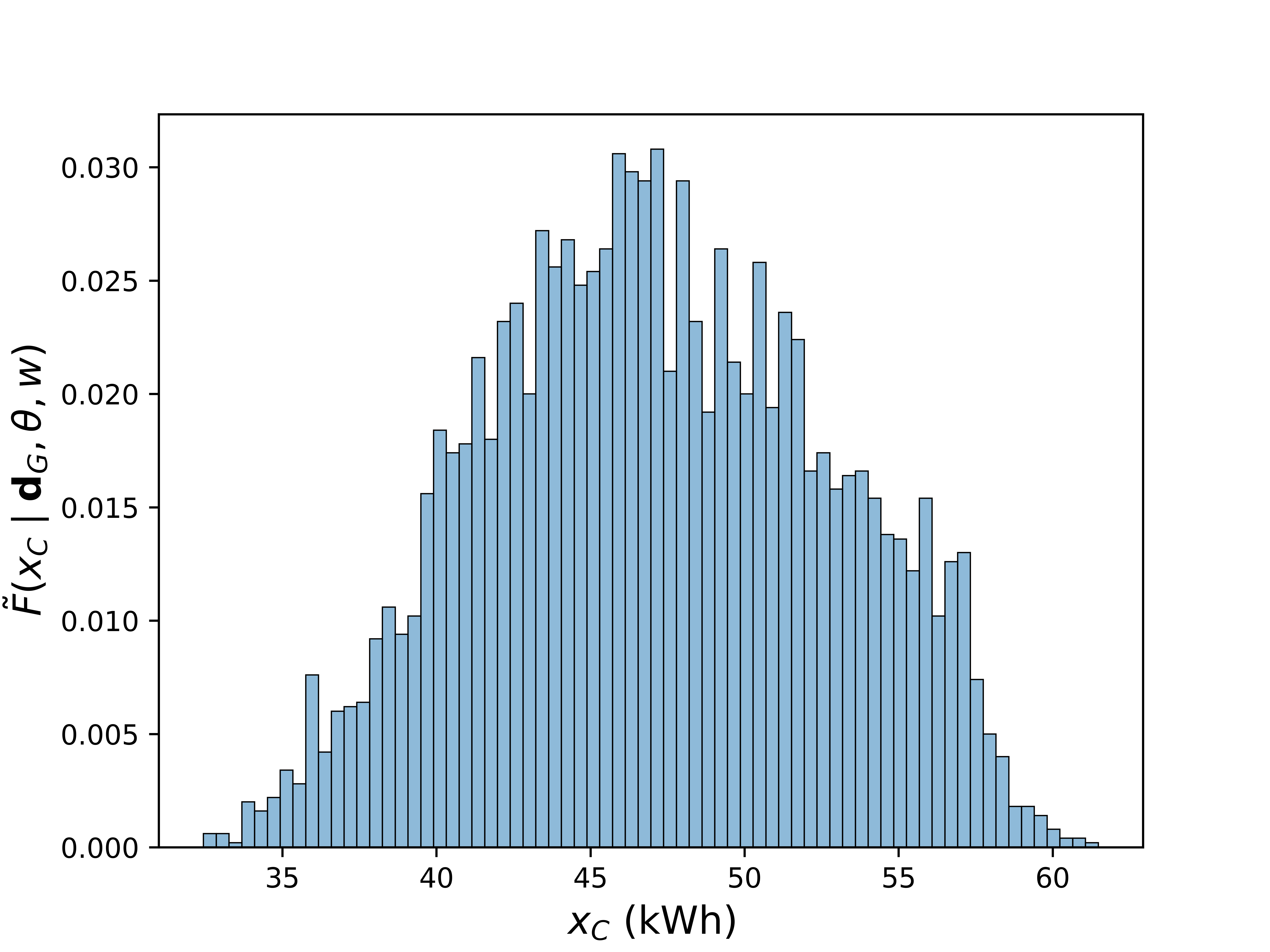}}
    \caption{Histograms  of  empirical   distributions  of  the   energy  drawn from the  battery, $\tilde{\mathsf{F}}(x_C \mid  {\mathbf d}_G, \theta, y)$, for (a) summer ($y=s$), and  (b) winter ($y=w$). 
    }
    \label{fig:histoXc}
\end{figure}

The physical model (\ref{First_Equation})-(\ref{eq:cumenergy}) relates $x_C$ {\em deterministically\/} to the {\em observed\/} data (i.e.\ the GPS record, ${\mathbf d}_G$),  the {\em known\/} EV-specific parameters, $\theta$ (Table~\ref{tab:vehicle_params}), the {\em known\/} season, $Y=y\in\{s,w\}$, but also to the {\em unknown\/} and unobserved quantities, $m_{peop}$ (\ref{eq:GMM}) and $W$ (\ref{eq:Gam}); i.e.
\[
x_C \equiv g({\mathbf d}_G, \theta, y, m_{peop}, W),
\]
with $g(\cdot)$ specified by (\ref{First_Equation})-(\ref{eq:cumenergy}). Hence, $x_C$ is unknown (i.e.\ a r.v.), whose distribution we now seek. Having elicited the marginal probability models of $m_{peop}$ and $W$, we assume---reasonably---that they are stochastically independent, i.e.\ $m_{peop} \indep W$. Formally, then, our entire knowledge base is processed consistently into {\em marginal\/} inference of unknown $x_C$---having marginalized over the unknowns---as follows:
\begin{multline}
\mathsf{F}(x_C \mid  {\mathbf d}_G, \theta, y)
= \label{eq:xCdistcont}\\[1.5mm]
\int_{0}^{\infty} \int_{0}^{\infty} \delta \left(x_C - g({\mathbf d}_G, \theta, y, m_{peop}, W)\right) \mathsf{F}(m_{peop})\mathsf{F}(W \mid y)dW dm_{peop}.
\end{multline}
We construct the standard empirical approximation~\cite{Tanner1996} of this distribution via a random sample of size $n \gg 0$, from $\mathsf{F}(m_{peop})$ (\ref{eq:GMM}) and $\mathsf{F}(W \mid y)$ (\ref{eq:Gam}):
\begin{equation}
\tilde{\mathsf{F}}(x_C \mid  {\mathbf d}_G, \theta, y)
\equiv
\frac{1}{n}\sum_{i}  \delta \left(x_C - g({\mathbf d}_G, \theta, y, m_{peop}^{(i)}, W^{(i)})\right), \label{eq:xCdist}
\end{equation}
where
\begin{eqnarray}
m_{peop}^{(i)} &\sim& \mathsf{F}(m_{peop}), \nonumber\\[1.5mm]
W^{(i)} &\sim& \mathsf{F}(W \mid y).\nonumber
\end{eqnarray}
We will use (\ref{eq:xCdist}) to predict $x_C$ for a particular EV model, $\theta$, with  GPS record, ${\mathbf d}_G$, operating in  a particular season, $Y=y\in\{s,w\}$. The seasonality of this predictive distribution of energy drawn from the EV battery is evident in the examples displayed in Figure~\ref{fig:histoXc}.  Note that the support of the empirical approximation (\ref{eq:xCdist}) is the random sample, $\{x_C^{(1)}, \ldots, x_C^{(n)}\}$, where
\begin{equation}
x_C^{(i)} \equiv g({\mathbf d}_G, \theta, y, m_{peop}^{(i)}, W^{(i)}), \;\; i=1, \ldots, n.
\label{eq:iidsample}
\end{equation}
In Fig.~\ref{fig:histoXc}, we display the    histograms formed by uniformly binning these iid (i.e.\ independent, identically distributed) samples to (quantized) values, $\hat{x}_C$. This reveals  the form of the underlying pdfs (\ref{eq:xCdistcont}). 









\subsection{Probability model of \(x_U\), the undeclared energy }
\label{sec:xU}

The net amount of undeclared energy, $x_U$, supplied to the EV battery during the certified  interval is, by definition, unknown. Its distribution has two components, conditional on the behaviour of the driver, as follows:
\begin{equation}
   \begin{array}{rr}
H_1: &\textsf{at least one undeclared charging event occurred}\;\; (\Pr[H_1] \equiv p_1),\\[4pt]
H_0\equiv \overline{H_1}: &\textsf{no undeclared charging event occurred}\;\; (\Pr[H_0] \equiv p_0 =1-p_1).
\end{array} 
 \label{eq:hypprior}
 \end{equation}
Here, $p_1$ is the {\em a priori\/} probability that the driver engaged in undeclared charging. The primary focus of this paper is to update this probability---using Bayesian inference, as explained in the sequel---to its data-driven  {\em a posteriori\/} value (\ref{eq:hyptest}), 
having processed the specific EV data from the last certified  interval.

In the $H_1$ case, we adopt the conservative uniform distribution in the interval, $I_U \equiv (0,\overline{x}_U]$ (kWh), where $\overline{x}_U$ is the maximum possible amount of undeclared energy\footnote{\label{foot:batt} For instance, $\overline{x}_U$ equals the  (energy) capacity of the EV battery, $E_{max}$ (Table~\ref{tab:vehicle_params}), if we assume that only one undeclared charging event occurred during the certified interval. We do not actively invoke this assumption.\\
On a more subtle point, the $x_U$ modelling in (\ref{eq:hypprior})-(\ref{eq:slabspike}) is elicited independently of the EV-specific database; i.e.\ $\mathsf{F}(x_U \mid {\mathbf d}_G, \theta, y) \equiv \mathsf{F}(x_U) $.
}:
\begin{equation}
   \mathsf{F}(x_U \mid H_1) \equiv {\mathcal U}_{I_U}(x_U). 
\label{eq:xUH1}
\end{equation}
Conversely, $x_U \stackrel{\mathsf a.s.}{=} 0$ in the $H_0$ case:
\begin{equation}
  \mathsf{F}(x_U \mid H_0) \equiv \delta(x_U).  
\label{eq:xUH0}
\end{equation}
Hence, the marginal model of $x_U$ is the slab-and-spike distribution:
\begin{equation}
\mathsf{F}(x_U) \equiv p_1{\mathcal U}_{I_U}(x_U) + (1-p_1) \delta(x_U).
\label{eq:slabspike}
\end{equation}

\subsection{Induced distribution of  certified differential SoC, $x_D$}

The energy conservation equation (\ref{eq:conser}) expresses the observed $x_D$ (i.e.\ the certified differential SoC) as the difference between the two unobserved energy quantities, i.e.\ $x_C \sim \tilde{\mathsf{F}}(x_C \mid \cdot)$ (\ref{eq:xCdist}) (the cumulative energy drawn from the battery) and $x_U \sim \mathsf{F}(x_U)$  (\ref{eq:slabspike}) (the undeclared energy), respectively. 
Analysis is greatly facilitated by our assumption of independence between these quantities. In this case, the induced {\em a posteriori\/} distribution of $x_D$---i.e.\ {\em after\/} processing the EV-specific evidence, ${\mathbf d}_G$, $\theta$ and $y\in\{s,w\}$, from the certified interval---is the {\em deterministic correlation\/} of $\tilde{\mathsf{F}}(x_C \mid \cdot)$ and $\mathsf{F}(x_U)$~\cite{AP91}:
\begin{eqnarray}
\tilde{\mathsf{F}}(x_D \mid {\mathbf d}_G, \theta, y) &=& \tilde{\mathsf{F}}(x_C \mid {\mathbf d}_G, \theta, y) \ast \mathsf{F}(-x_U)|_{x_C \rightarrow x_D, x_U \rightarrow x_D}\nonumber\\[1.5mm]
&=& p_1 \tilde{\mathsf{F}}(x_C \mid {\mathbf d}_G, \theta, y)|_{x_C \rightarrow x_D}\ast {\mathcal U}_{I_U} (-x_D) + (1-p_1)  \tilde{\mathsf{F}}(x_C \mid {\mathbf d}_G, \theta, y)|_{x_C \rightarrow x_D} \nonumber \\[1.5mm]
&=& p_1 \tilde{\mathsf{F}}(x_D \mid {\mathbf d}_G, \theta, y, H_1)+(1-p_1)\tilde{\mathsf{F}}(x_D \mid {\mathbf d}_G, \theta, y, H_0).\label{eq:xD2comp}
\end{eqnarray}
This is a binary mixture model with components,
\begin{eqnarray}
\tilde{\mathsf{F}}(x_D \mid {\mathbf d}_G, \theta, y, H_1) &\equiv& \tilde{\mathsf{F}}(x_C \mid {\mathbf d}_G, \theta, y)|_{x_C \rightarrow x_D}\ast {\mathcal U}_{I_U} (-x_D),\label{eq:xD1stcomp}\\[1.5mm]
\tilde{\mathsf{F}}(x_D \mid {\mathbf d}_G, \theta, y, H_0) &\equiv& \tilde{\mathsf{F}}(x_C \mid {\mathbf d}_G, \theta, y)|_{x_C \rightarrow x_D}. \label{eq:xDH0}
\end{eqnarray}
These components are conditioned on the hypothesis of interest (\ref{eq:hypprior}) being true, i.e.\ $H_1$, and false, i.e.\ $H_0$, respectively. The mixing probabilities, $p_1$ and $1-p_1$, are the {\em a priori\/} probabilities of these hypotheses (see footnote~\ref{foot:batt}), and $\ast$ denotes linear convolution. For convenience, we compute $\tilde{\mathsf{F}}(x_D \mid {\mathbf d}_G, \theta, y, H_1)$  by discretizing ${\mathcal U}_{I_U} (x_U)$---denoted by $\tilde{\mathcal U}_{I_U} (x_U)$---to the same bins, $\hat{x}_C$---as those used for $\tilde{\mathsf{F}}(x_C \mid \cdot )$ (\ref{eq:xCdist}) (see Fig.~\ref{fig:histoXc}), yielding the discrete correlation,
\begin{equation}
\tilde{\mathsf{F}}(\hat{x}_D \mid \cdot, H_1) = \tilde{\mathsf{F}}(\hat{x}_D \mid \cdot )\ast {\mathcal U}_{I_U} (-\hat{x}_D) \equiv \sum_{\hat{x}_U \in I_U} 
\tilde{\mathsf{F}}(\hat{x}_U \mid \cdot )\times \tilde{\mathcal U}_{I_U} (\hat{x}_U + \hat{x}_D ).
\label{eq:disconv}
\end{equation}

\begin{figure}[H]
    \centering
    \subfigure[]{\includegraphics[width=0.45\linewidth]{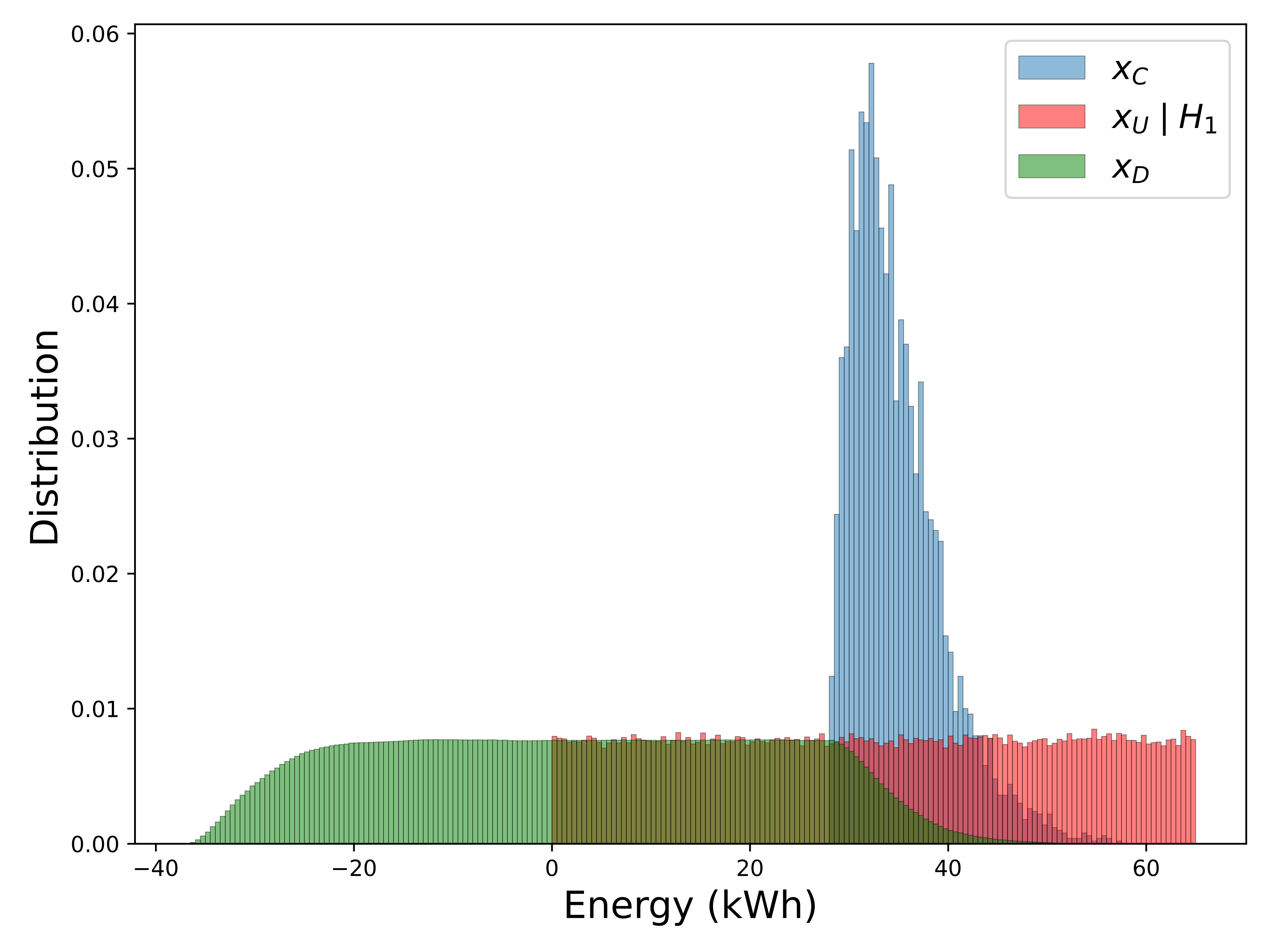}}
    \subfigure[]{\includegraphics[width=0.45\linewidth]{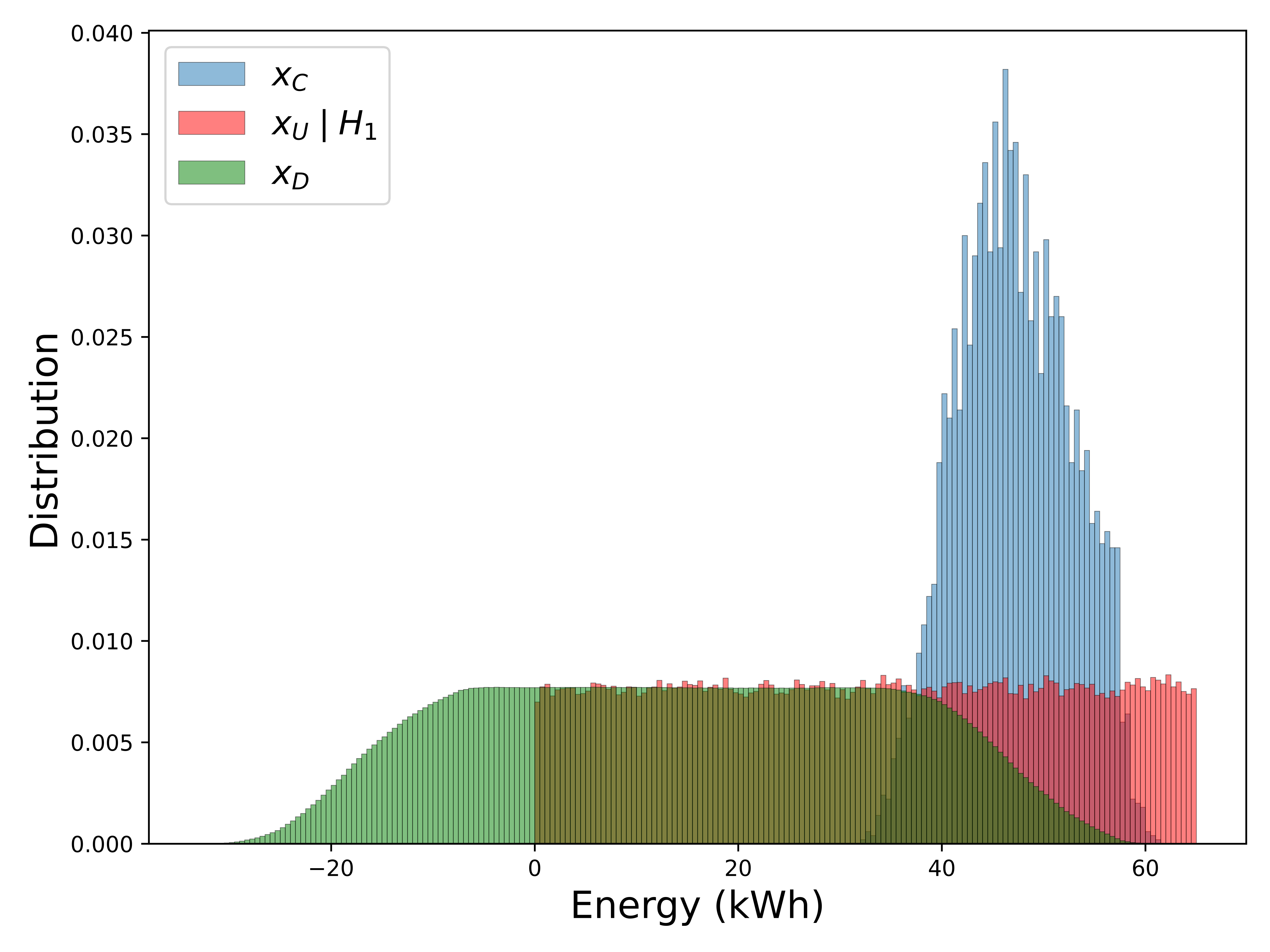}}
    \caption{Histogram representations  (in green) of  the empirical   distributions, $\tilde{\mathsf{F}}(x_D \mid  {\mathbf d}_G, \theta, y)$ \eqref{eq:xD2comp},  of  the certified differential SoC, $x_D$, in (a) summer ($y=s$), and  (b) winter ($y=w$). Also displayed are the contributing distributions in \eqref{eq:xD2comp}, i.e.\ $\tilde{\mathsf{F}}(x_C \mid  {\mathbf d}_G, \theta, y)$ (\ref{eq:xCdist}) (in blue) and the conditional distribution, $\tilde{\mathsf{F}}(x_U \mid  H_1)$ (in red).}
    \label{fig:Convolution}
\end{figure}

\begin{remark}[Range of $x_D$] Let the actual range of $x_C >0$---i.e. the support of distribution (\ref{eq:xCdist})---be $I_{C} \equiv (\underline{x}_C, \overline{x}_C)$, and recall that $x_U \in I_U \equiv (0, \overline{x}_U)$. It follows from (\ref{eq:xD2comp}) that $x_D \in I_{D} \equiv (\underline{x}_D, \overline{x}_D) = (\underline{x}_C-\overline{x}_U,\;  +\overline{x}_C)$. This is best understood by noting that the $H_1$-component  in   (\ref{eq:xD2comp})---i.e.\ (\ref{eq:xD1stcomp}), implemented via (\ref{eq:disconv})---implements a sliding sum of $\tilde{\mathsf{F}}(\hat{x}_U \mid \cdot)$. These bounds for $x_D$ capture the extreme cases reasonably: (i) if $x_C =0$ (\ref{eq:conser}), then  $\overline{x}_U = E_{max}$ (see footnote~\ref{foot:batt}), and so $\underline{x}_D = -x_{U}=-E_{max}$. Conversely, if $x_U=0$, then $\overline{x}_C = E_{max}$, and so $\overline{x}_D = \overline{x}_C = E_{max}$. 
\end{remark}
Figure~\ref{fig:Convolution} displays the binned empirical distributions (i.e.\ histograms) of $x_D$, corresponding to the cases in Figure~\ref{fig:histoXc}.  

\subsection{Testing the hypothesis of undeclared charging, $H_1$} 
We now compute a Bayesian hypothesis test for $H_1$ (Section~\ref{sec:xU}) via the  binary mixture model of $x_D$ (\ref{eq:xD2comp}). As explained near the beginning of Section~\ref{sec:hypo}, this involves the sequential\footnote{The inference is sequential, in that have already processed  the EV-specific database from the certified interval, i.e.\ $({\mathbf d}_G, \theta, y)$. In this sense, $\tilde{\mathsf{F}}(x_D \mid {\mathbf d}_G, \theta, y)$ (\ref{eq:xD2comp}) is a data-informed prior for $X_D$, and we wish then to soft-classify (i.e.\ infer the label of) the observed $X_D = x_D$.} processing of (i.e.\ conditioning on) the measured differential SoC, $x_D$ (\ref{eq:conser}), yielding  $\Pr[H_1 \mid {\mathbf d}_G, \theta, y, x_D]$ (\ref{eq:hyptest}).  
It follows from a simple application of Bayes' rule~\cite{AP91}:
\begin{equation}
    \label{eq:Bayes}
\Pr[H_1 \mid {\mathbf d}_G, \theta, y, x_D] = \frac{p_1 \tilde{\mathsf{F}}(x_D \mid {\mathbf d}_G, \theta, y, H_1)}
{p_1 \tilde{\mathsf{F}}(x_D \mid {\mathbf d}_G, \theta, y, H_1) + (1-p_1) \tilde{\mathsf{F}}(x_D \mid {\mathbf d}_G, \theta, y, H_0)}.
\end{equation}
The components on the right-hand-side of (\ref{eq:Bayes}) are evaluated at the observed $x_D$ via (\ref{eq:xD1stcomp}) and (\ref{eq:disconv}).   








\section{Simulation results for the Bayesian test of the  undeclared charging hypothesis, $H_1$}
\label{sec:sims}

In this section, we present  extensive simulations  to validate the Bayesian hypothesis testing  framework---which we will call the ``algorithm''---defined via (\ref{eq:Bayes}).   For this purpose, we deploy the popular mobility simulator, SUMO, to reproduce a road network and its traffic. We use the resulting simulated data to validate the algorithm  under controlled conditions. Section \ref{SUMO}  briefly presents the mobility simulator and the selected case study, while Section \ref{Simulation_Results} reports the performance of the algorithm under different operating conditions.

\subsection{The SUMO mobility simulator }
\label{SUMO}

SUMO (Simulation of Urban MObility)~\cite{kurczveil2014implementation} is a microscopic and continuous multi-modal traffic simulation platform that has been designed to handle and analyse large traffic networks. Several key inputs are required, which we instantiate in our current work as described next: 

\begin{itemize}
    \item {\em Simulation environment:} our simulated road network covers  part of \textit{central London}. It is imported using OpenStreetMap \cite{OpenStreetMap}, an editable map database. We have divided the road network into 11 distinct areas, as depicted in Figure \ref{fig:London}.
    
    \begin{figure}[H]
    \centering
    \includegraphics[width=0.8\linewidth]{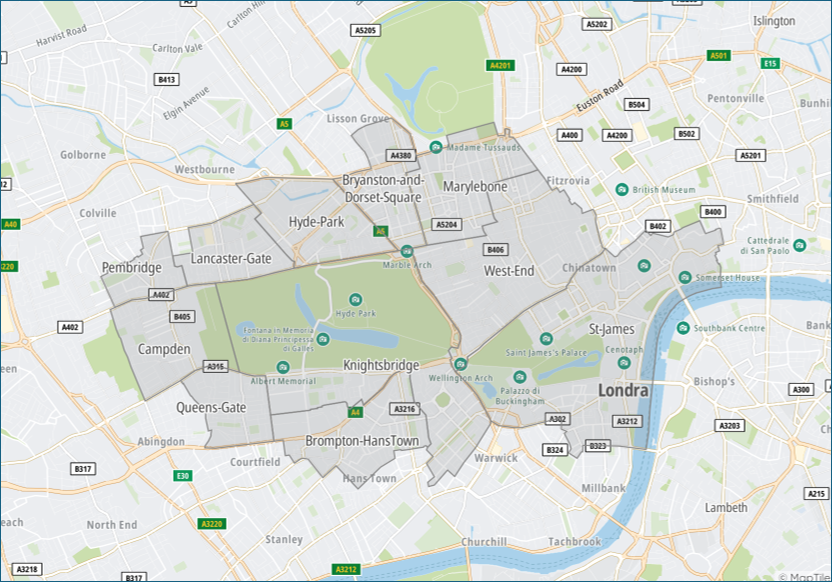}
    \caption{Map of the neighbourhoods of London (UK) which have been used for our case study. The map has been imported through OpenStreetMap~\cite{OpenStreetMap}.}
    \label{fig:London}
    \end{figure}

    \item {\em EV properties:} recall that the algorithm must have access to the parameters, $\theta$,  of the specific EV (Table~\ref{tab:vehicle_params}). For simplicity, we adopt $\theta$ for the KIA Soul EV 2020 (Table \ref{tab:vehicle_params_specifc}), this being the default  compact EV model  in the {\em electric\/} class of SUMO.  

\begin{table}[H]
\centering
\begin{tabular}[t]{r|c|l}
\textbf{Parameter} & \textbf{Symbol} & \textbf{Value } \\ \hline
Battery capacity & $E_{max}$  & 35 kWh  \\ \hline
{Kerb mass} & $m_{veh}$ & 1682 kg  \\ \hline
Front surface area & $A_{veh}$ & 2.6 m$^2$  \\ \hline
Internal moment of inertia & $J_{int}$ & 40 
kg.m$^2$ \\ \hline
Radial drag coefficient & $c_{rad}$ & 0.1  \\ \hline
Roll drag coefficient & $c_{roll}$ &   0.01 \\ \hline
Air drag coefficient & $c_{w}$ &   0.35 \\ \hline
Propulsion efficiency & $\eta_{prop}$ & 0.98   \\ \hline
Recuperation efficiency & $\eta_{recup}$ & 0.96   \\ 
\end{tabular}
\caption{Specific EV parameters (for the KIA Soul EV 2020), $\theta$, used in the  SUMO-based simulations.}
\label{tab:vehicle_params_specifc}
\end{table}

    \item {\em Origin-destination (OD) matrix:} the entries in this $11\times 11$ matrix record the numbers of trips between the 11 origin areas---indexing the rows of the matrix---and (the same 11) destination areas---i.e.\ matrix columns---in central London (Figure~\ref{fig:London}). The OD matrix is used to simulate traffic flows by specifying the number of EVs that travel between the areas of the  road network in a specific time frame. An example of this OD matrix for one hour is provided in Figure \ref{fig:OD_Matrix}. In principle, these OD matrices are available to city municipalities via sensing infrastructure and surveys. The OD matrix in Figure~\ref{fig:OD_Matrix} was generated using the TomTom O/D Analysis API \cite{TomTomODAnalysis}. This API analyzes historical traffic data using\textit{ floating car data} (FCD) from millions of anonymized TomTom navigation devices and apps, providing realistic insights into traffic flows and patterns. The numbers of trips reported in Figure \ref{fig:OD_Matrix} have been proportionally down-scaled to simulate traffic-free conditions.
          
    \begin{figure}[H]
    \centering
    \includegraphics[width=0.9\linewidth]{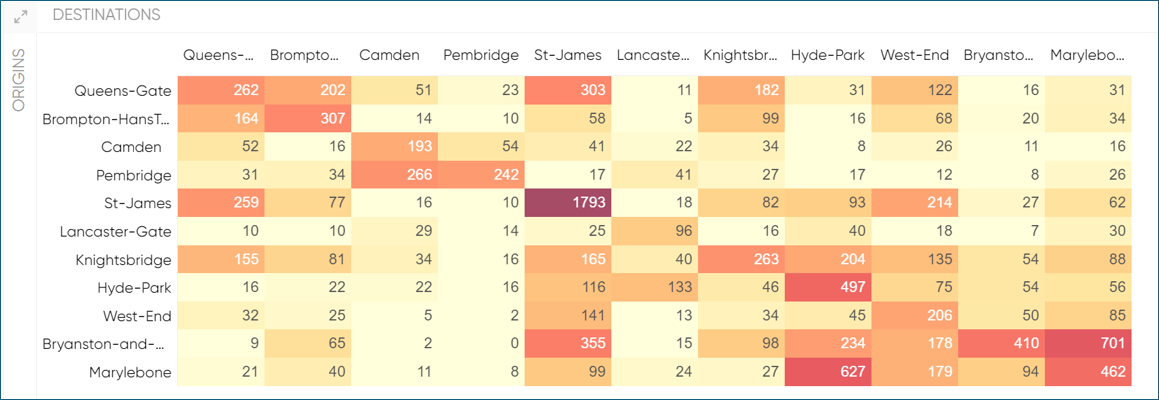}
    \caption{The OD matrix generated via the {\em TomTom O/D Analysis API}~\cite{TomTomODAnalysis} for the time frame 08:00-09:00 on Monday 24th October 2022. }
    \label{fig:OD_Matrix}
    \end{figure}
\end{itemize}

Recall, from Section~\ref{sec:physmodel}, that the cumulative energy drawn from the battery during a trip of $t$ s,  i.e.\ $x_C[t]$ (\ref{eq:cumenergy}), is evaluated in discrete time, $t\in\{0,1,2,\ldots\}$ s, via the physical model (\ref{First_Equation})-(\ref{Final_Equation}). This is consistent with the model deployed by the \textit{electric} class in SUMO. 
Hence, the SUMO simulations are conducted in discrete time, using a one-second time step (Table~\ref{tab:GPS_dataframe}). 






\subsection{Performance of the hypothesis test for undeclared charging }
\label{Simulation_Results}

\begin{figure}[H]
    \centering   \includegraphics[width=0.7\linewidth]{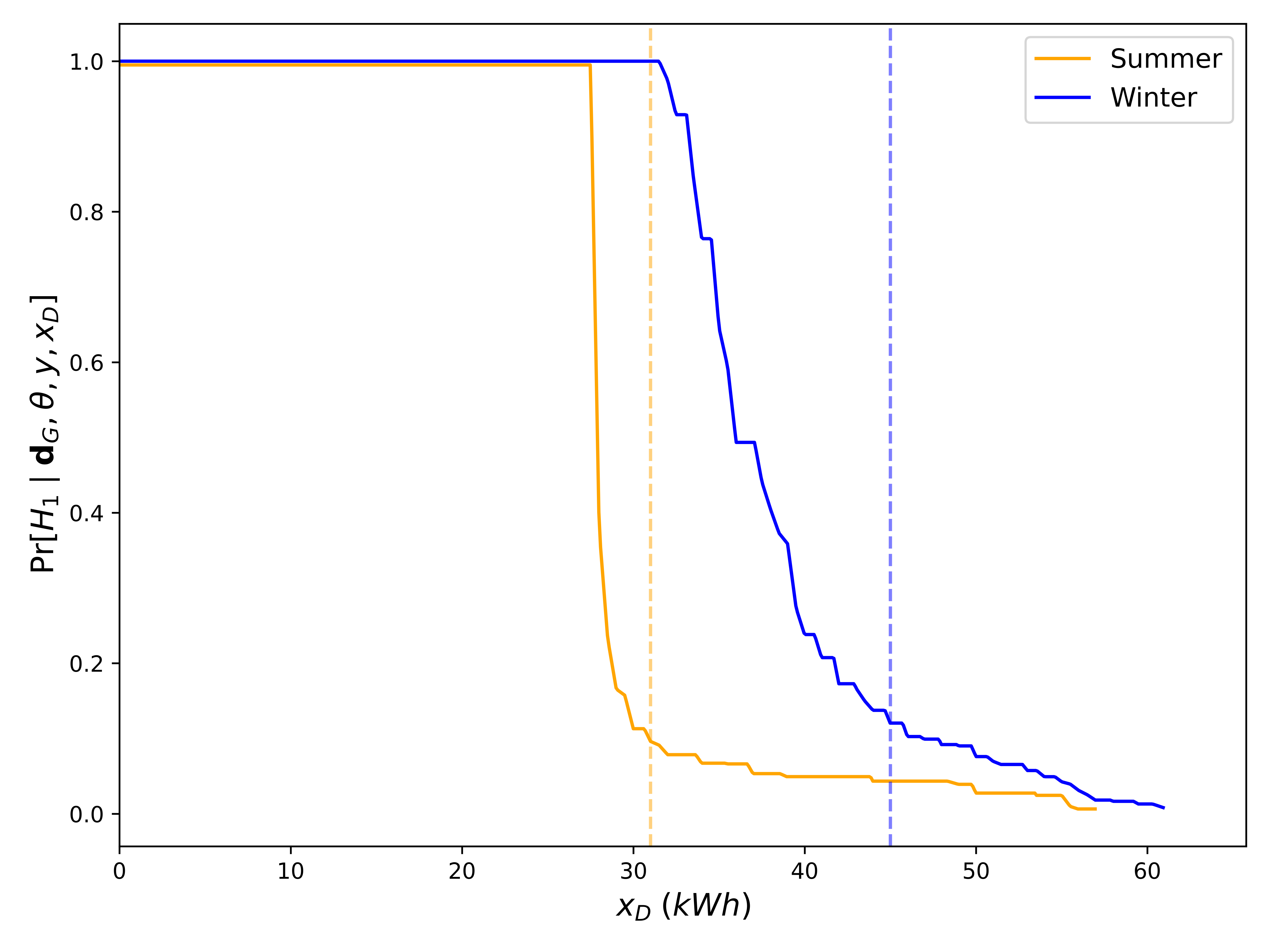}
    \caption{Probability of undeclared charging, $\Pr[H_1 |\cdot, x_D]$,  as a function of the certified differential SoC, $x_D$. The dotted lines indicate $x_D \equiv x_{C,\mathsf{MAP}}$ for the respective seasons (see Figure~\ref{fig:histoXc}).} 
    \label{Hypothesis_Testing}
\end{figure}

In Figure \ref{Hypothesis_Testing}, we display $\Pr[H_1 |{\mathbf d}_G, \theta, y, x_D] \equiv \Pr[H_1 |\cdot, x_D]$ (\ref{eq:Bayes}) (being the probability that at least one undeclared charging event occurred during the certified interval) as a function of $x_D$ (\ref{eq:conser}) (i.e.\ the measured differential SoC for the  interval). The certified interval duration, $T_C$ (Section~\ref{sec:hypo}), is   set at about  two weeks, corresponding to about 40 trips, each about 8 km long (i.e.\ a total of about 320 km in two weeks). 

The dependence of $\Pr[H_1 |\cdot, x_D]$, on the season, $Y=y\in\{s,w\}$, is evident in Figure~\ref{Hypothesis_Testing}. In both cases, $\Pr[H_1 |\cdot, x_D]\rightarrow 1$,   for  low values of $x_D$, i.e.\ for those differential SoCs which are inconsistent with the  GPS record (Table~\ref{tab:GPS_dataframe}), and, therefore, inconsistent with the values of $x_C$ that are probable under (\ref{eq:xCdist}).    Meanwhile, $\Pr[H_1 |\cdot, x_D]\rightarrow 0$  as $x_D$ increases to values with high probability under $H_0$ (\ref{eq:xDH0}). Note the more rapid transition in $\Pr[H_1 |\cdot, x_D]$ in summer than in winter. This is caused by the increased uncertainty in $x_C$ (\ref{eq:xCdist}) in winter (see Figure~\ref{fig:histoXc}), which, in turn, is induced by the increased uncertainty of $W$ in  winter (i.e.\ in the  power drawn by the heating system) compared to $W$ in summer (i.e.\ drawn by the A/C system). In Figure~\ref{Hypothesis_Testing}, we have indicated the seasonal cases when $x_D = x_{C,\mathsf{MAP}}$, i.e. the modal values of $x_C$, respectively (see Figure~\ref{fig:histoXc}). This verifies the satisfactory performance of the algorithm, which assigns  high probabilities  to $H_0$ (i.e.\ to the hypothesis that no undeclared charging occurred) in this case. An empirical  study of the performance of the algorithm follows in Section~\ref{sec:MC}.

\subsection{The hypothesis test in representative scenarios }
We now describe how $\Pr[H_1 |\cdot, x_D]$ (\ref{eq:Bayes})---which has been graphed as a function of the measured differential SoC, $x_D,$ for the two seasonal cases, $y\in\{s,w\}$, in Figure \ref{Hypothesis_Testing}---can be deployed  in  the detection of various non-compliant behaviours (i.e.\ uncertified charging) by drivers. 

\begin{figure}[H]
    \centering
    \subfigure[]{\includegraphics[width=0.45\linewidth]{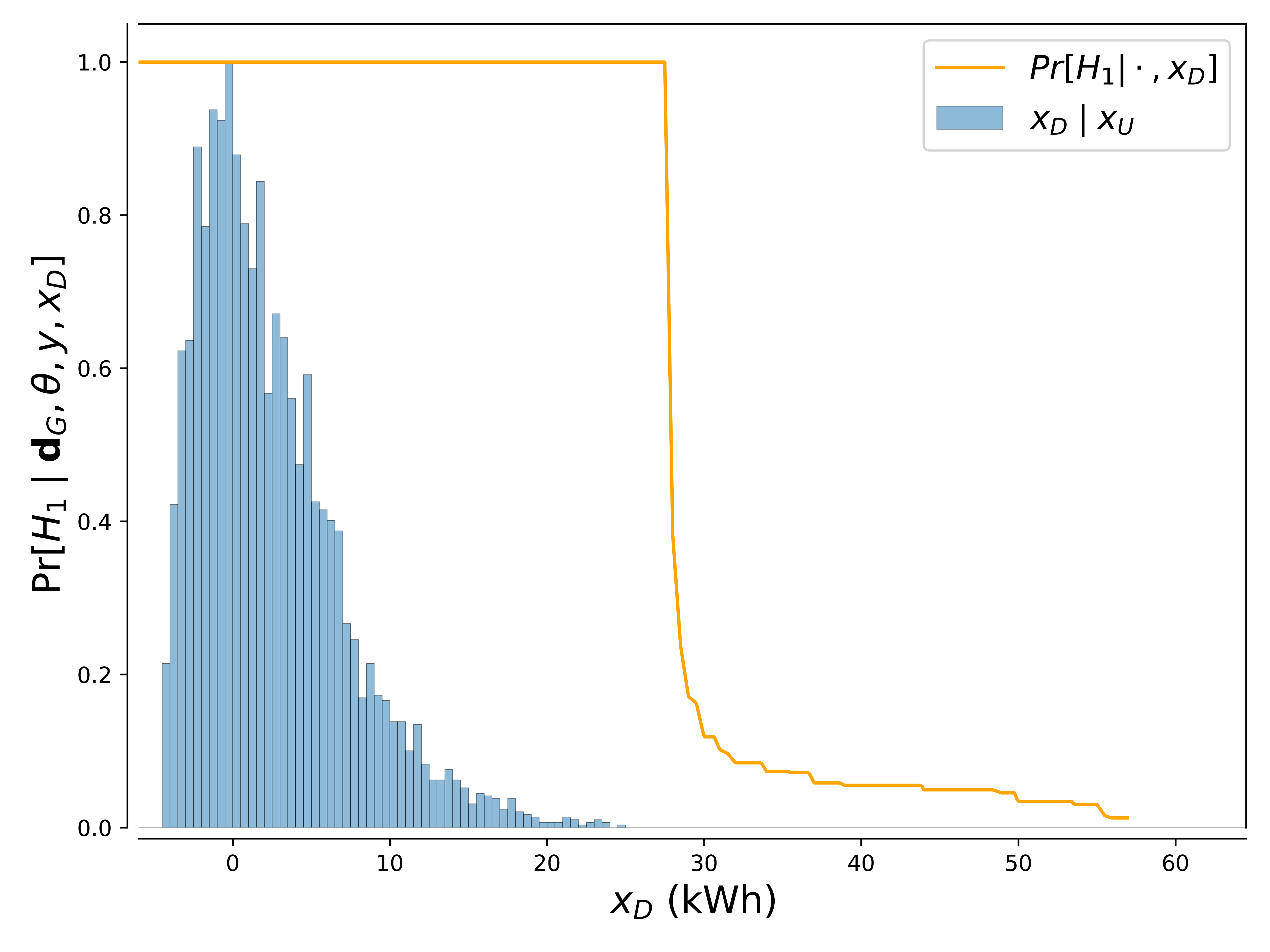}}
    \subfigure[]{\includegraphics[width=0.45\linewidth]{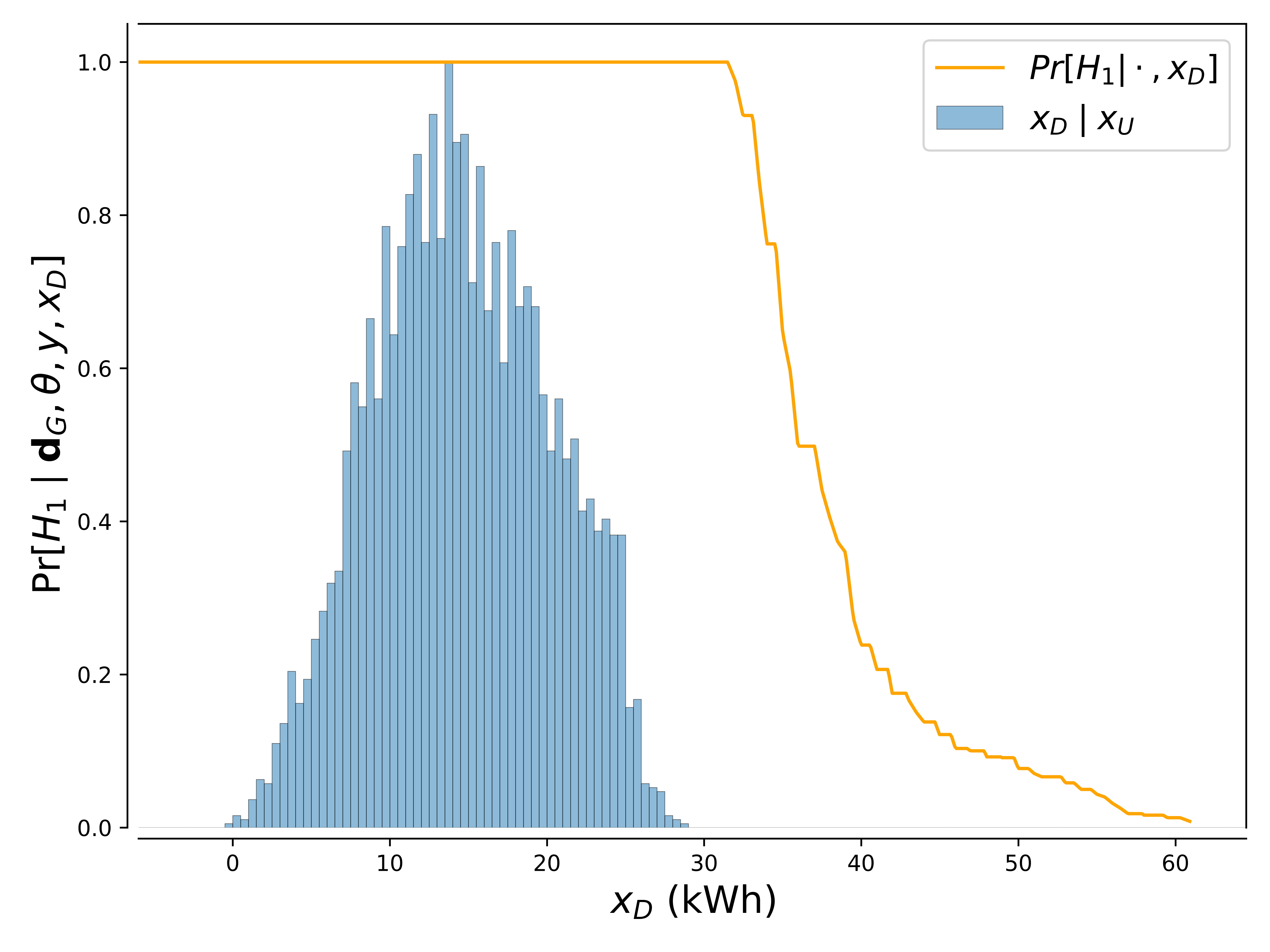}}
    \caption{Probability of undeclared charging, $\Pr[H_1 |\cdot, x_D]$, when {\em significant undeclared charging\/} has occurred (i.e.\ $H_1$, with $x_U= 0.5\,E_{max}$), in (a) summer, and (b) winter. The empirical distribution of $x_D\mid x_U$   is superimposed in each case.  } 
    \label{fig:PerformanceCheating50}
\end{figure}

\subsubsection{Scenario~1: $H_1$ true, with $x_U = 0.5\,E_{max}$ (Figure~\ref{fig:PerformanceCheating50})}
Here, we simulate $x_D$ when undeclared charging has, indeed, taken place, the amount of undeclared charge being 50\% of the battery capacity, $E_{max}$ (Table~\ref{tab:vehicle_params_specifc}: see Footnote~\ref{foot:batt}). In this flagrant non-compliance case, the measured differential SoC \eqref{eq:conser} is probably  small. Indeed, there is a strictly positive probability that $x_D <0$, as seen in summer (Figure~\ref{fig:PerformanceCheating50}(a)). Such values of $x_D$ provide trivial tests for the algorithm (see Footnote~\ref{foot:trivial}). Meanwhile, the hypothesis test correctly detects the undeclared charging event(s) for all simulated cases of $x_D$; i.e.\ note that $\Pr[H_1 |\cdot, x_D]\to 1$ across the entire range of simulated $x_D$ in this scenario.   In Section~\ref{sec:MC}, we will quantify the sensitivity of the test.


\begin{figure}[H]
    \centering
    \subfigure[]{\includegraphics[width=0.45\linewidth]{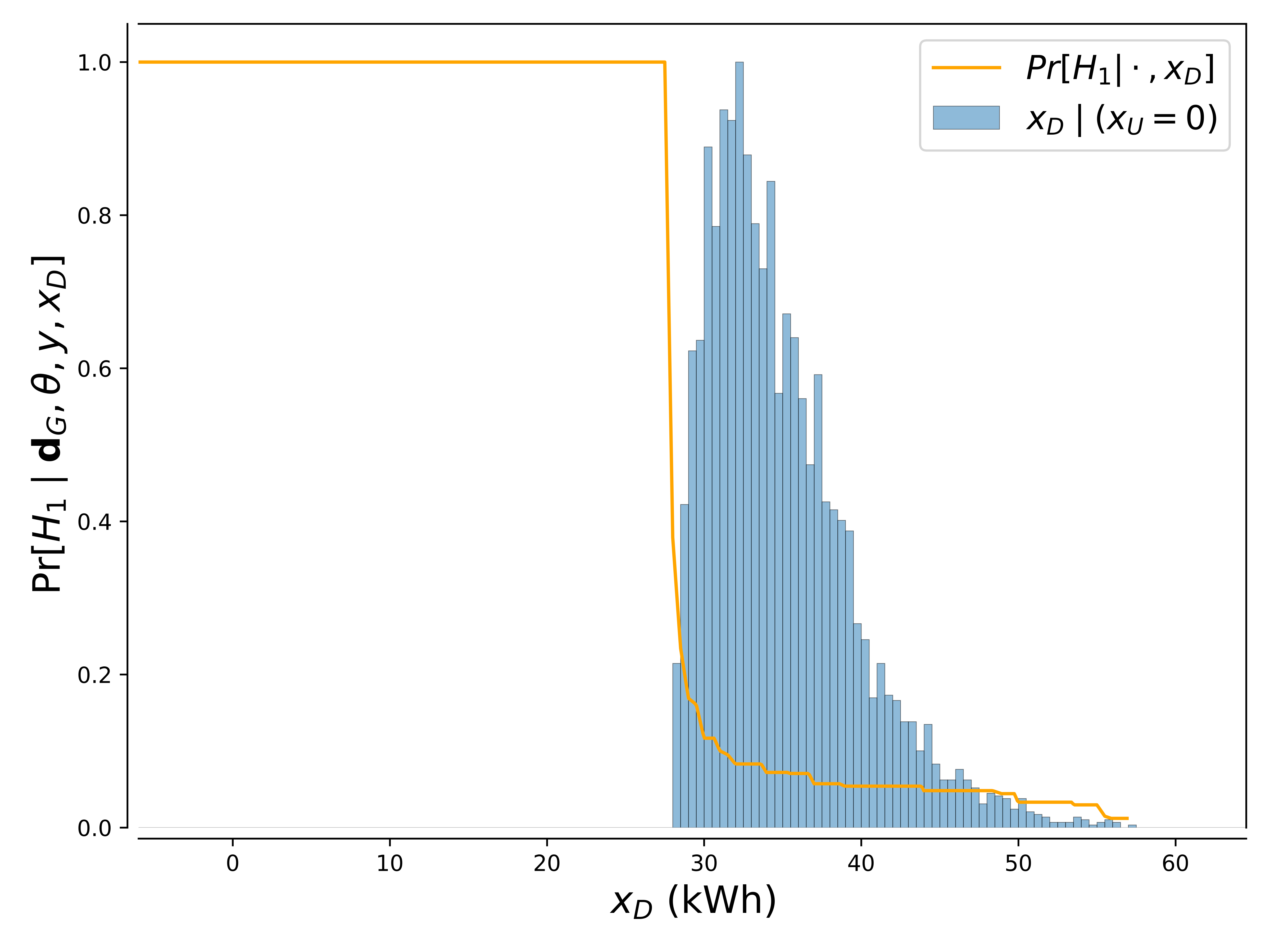}}
    \subfigure[]{\includegraphics[width=0.45\linewidth]{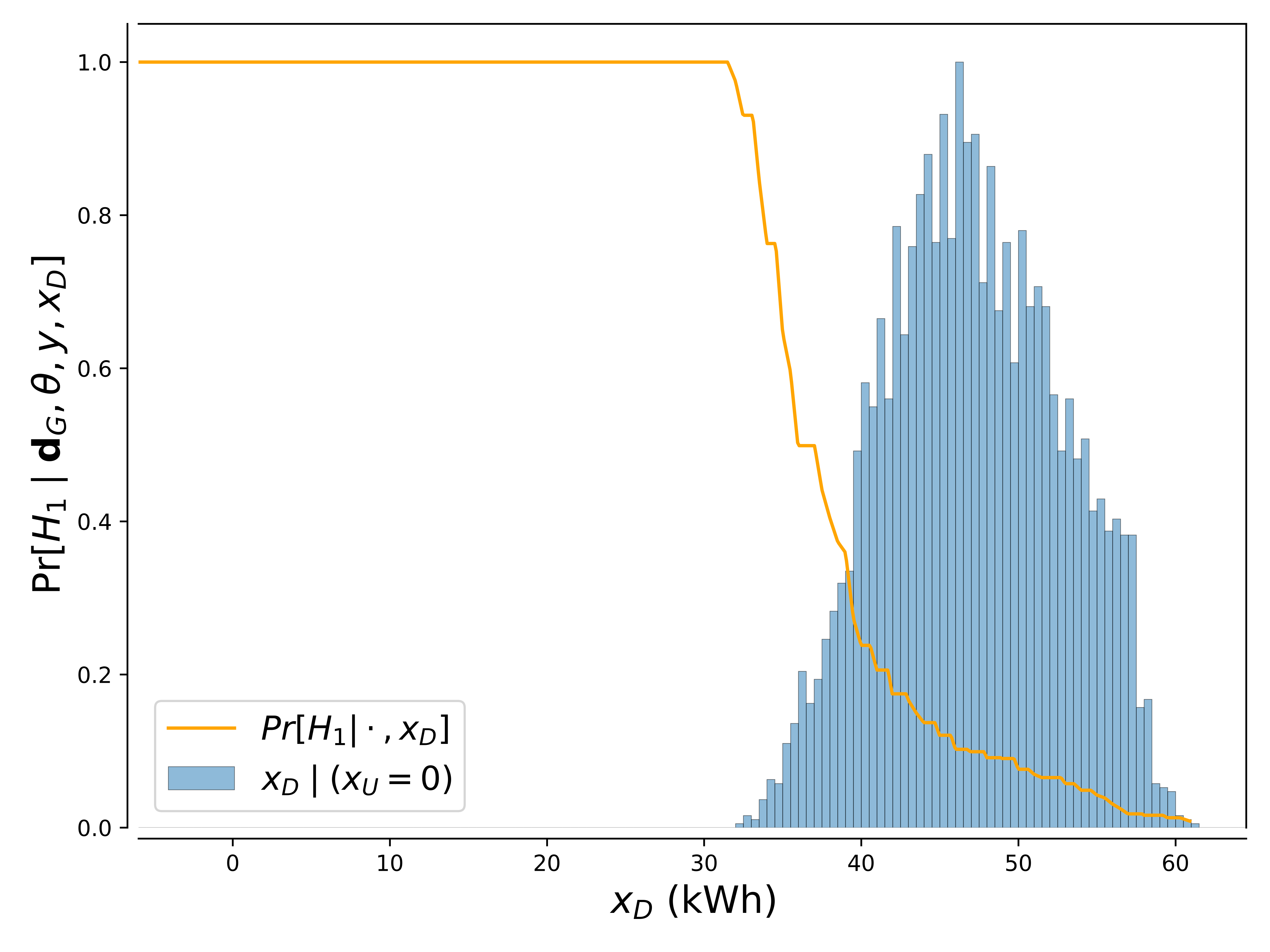}}
    \caption{Probability of undeclared charging, $\Pr[H_1 |\cdot, x_D]$, when there has been {\em no undeclared charging\/} (i.e.\ $H_0$ $\equiv$ $(x_U= 0)$),  in (a) summer, and (b) winter. The empirical distribution of $x_D\mid (x_U=0)$   is superimposed in each case. They are respectively equal to $\tilde{\mathsf{F}}(x_C \mid  {\mathbf d}_G, \theta, y)$ \eqref{eq:xCdist} in this $H_0$ case. }  
    \label{fig:PerformanceNoCheating}
\end{figure}

\subsubsection{Scenario~2: $H_0$ true, i.e.\ $x_U = 0 $ (Figure~\ref{fig:PerformanceNoCheating})}
We repeat the investigation when the driver has definitely been compliant with the charge certification scheme, in the sense that they have not charged their EV during the certified interval. The same GPS data, ${\mathbf d}_G$, are adopted as in scenario~1. Hence, the consumed energy distribution,  $\tilde{\mathsf{F}}(x_C \mid  {\mathbf d}_G, \theta, y)$ (\ref{eq:xCdist}), is the same as in scenario~1. For the same reason, the Bayesian hypothesis test, $\Pr[H_1 |{\mathbf d}_G, \theta, y, x_D]$ (\ref{eq:Bayes}) is also invariant. The test performs excellently once again; i.e.\ $\Pr[H_1 |\cdot, x_D]\to 0$ for most simulated cases of $x_D$ when $x_U =0$. This is so because  these $x_D$s  are, by definition (\ref{eq:conser}), iid draws from the $x_C$ predictor.  Nevertheless, false positives do occur, particularly in winter. The core reason for this is the increased uncertainty around the use of auxiliary services, i.e.\ $W$, evident in Figure~\ref{fig:Gam}. There is a significant probability of low $W$ \eqref{eq:DelEaux} in winter (i.e.\ drivers may choose not to use the heating system during mild periods), and, therefore, low $x_D$, in which case $\Pr[H_1 |\cdot, x_D] \gg 0$. We will assess the specificity of the test in  Section~\ref{sec:MC}.


\begin{figure}[H]
    \centering
    \subfigure[]{\includegraphics[width=0.45\linewidth]{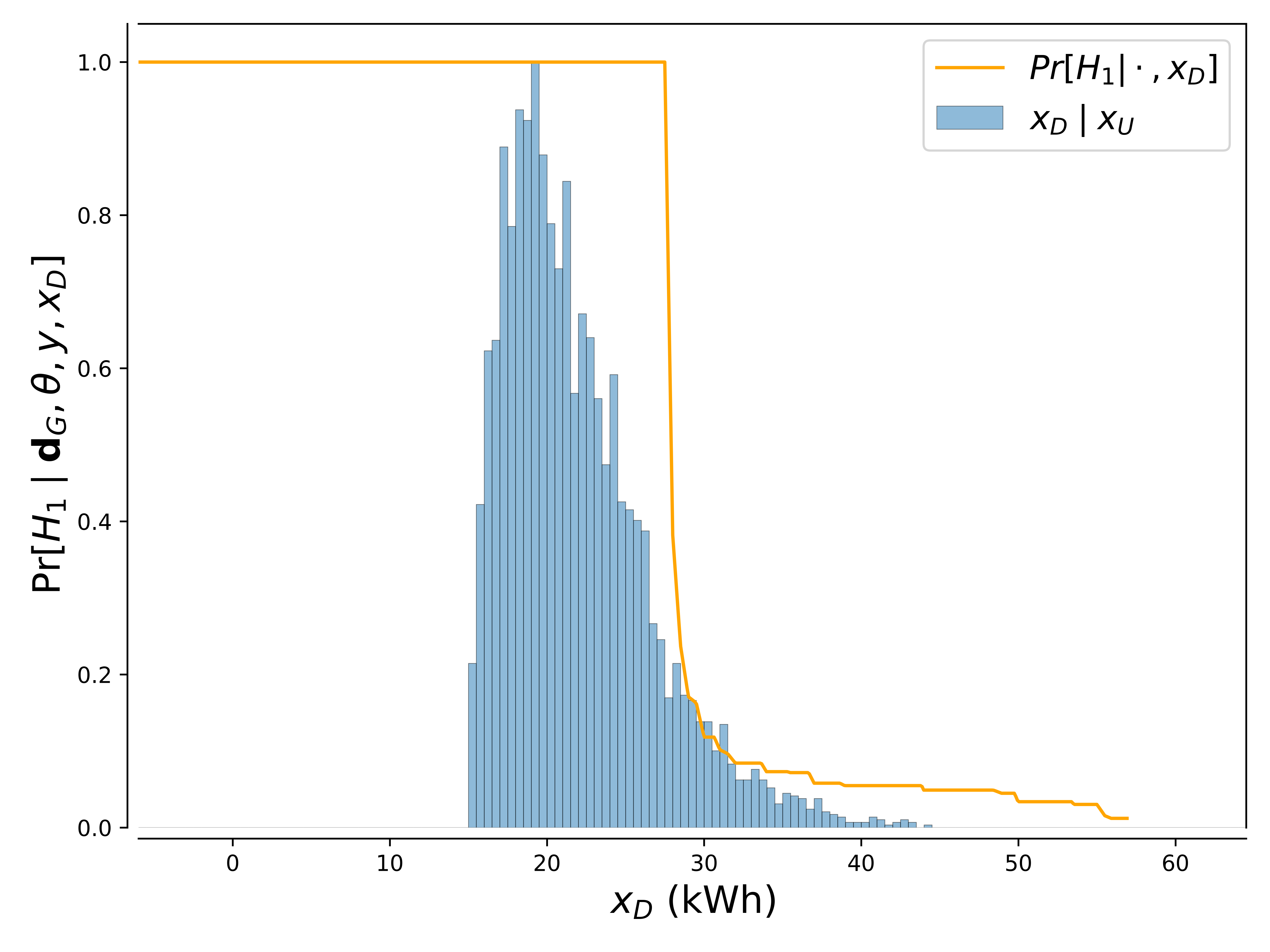}}
    \subfigure[]{\includegraphics[width=0.45\linewidth]{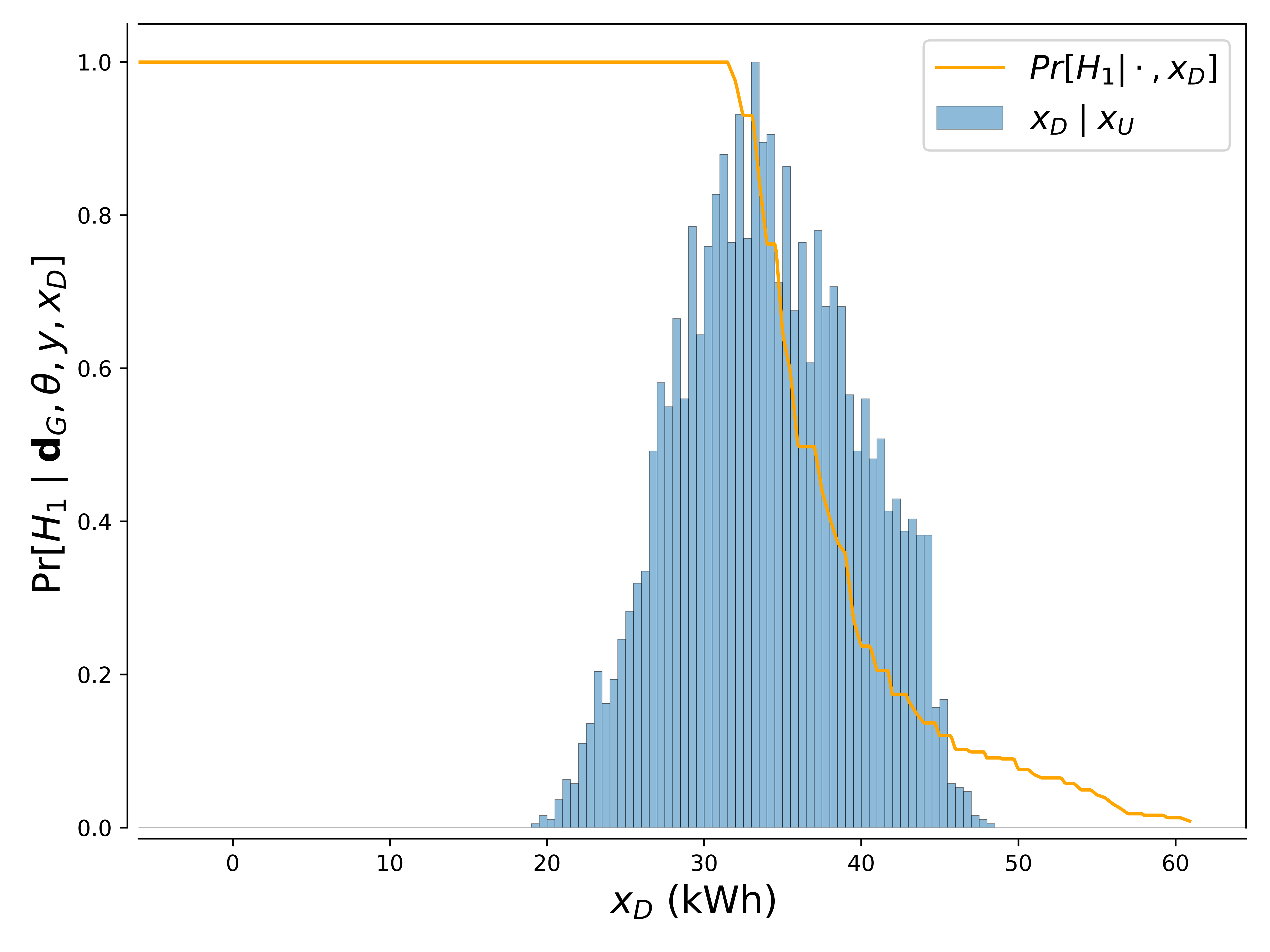}}
    \caption{Probability of undeclared charging, $\Pr[H_1 |\cdot, x_D]$, when {\em minor undeclared charging\/} has occurred (i.e.\ $H_1$, with $x_U= 0.2\,E_{max}$), in (a) summer, and (b) winter. The empirical distribution of $x_D\mid x_U$   is superimposed in each case.   }
    \label{fig:PerformanceCheating20}
\end{figure}

\subsubsection{Scenario~3:  $H_1$ true, with $x_U = 0.2\,E_{max}$ (Figure~\ref{fig:PerformanceCheating20})}
\label{sec:scen3}

Finally, we illustrate in Figure \ref{fig:PerformanceCheating20} the situation when drivers engage in a small amount of undeclared charging, i.e.  $x_U = 0.2\,E_{max}$.  Once again, the GPS data, ${\mathbf d}_G$, remain invariant. 
These---perhaps deliberate---cases of minor non-compliance with the certification scheme are unsurprisingly harder to detect, particularly in winter, when uncertainty around the use of auxiliary devices causes  $\Pr[H_1 |\cdot, x_D] \ll 1$ in a significant proportion of the simulated cases of $x_D$. 

In summary, the behaviour of the Bayesian hypothesis test is evident from these scenarios, all of which involve the same SUMO-simulated EV trips during the certified interval. The empirical distribution of $x_D$ is invariantly shifted by the amount of $x_U \geq 0$; i.e.\ inserting the energy conservation equation \eqref{eq:conser} into \eqref{eq:xD2comp}, then
\[
\tilde{\mathsf{F}}(x_D \mid {\mathbf d}_G, \theta, y, x_U) =  
\left.\tilde{\mathsf{F}}(x_C \mid {\mathbf d}_G, \theta, y)\right|_{x_C \rightarrow x_D+x_U}.
\]
Meanwhile, the Bayesian classifier,  $\Pr[H_1 |{\mathbf d}_G, \theta, y, x_D]$ (\ref{eq:Bayes}) is, of course,  not a function of $x_U$, since the latter has been marginalized out via \eqref{eq:xD1stcomp} and \eqref{eq:xDH0}. 


\subsection{Statistical Performance of the Test for Undeclared Charging }
\label{sec:MC}
Having gained an understanding of the behaviour of the testing algorithm for specific amounts of undeclared charging ($x_U >0$ given $H_1$  and $x_U=0$ given $H_0$ (\ref{eq:hypprior})), we want to assess its statistical performance---such as its sensitivity and specificity---under predictable behaviours of the driver. For this purpose, we hierarchically draw iid samples from the slab-and-spike mixture model of $x_U$ (\ref{eq:slabspike}) in a Monte Carlo (MC) simulation. We first simulate from the noninformative (i.e.\  unbiased) Bernoulli prior for $H_1$ (\ref{eq:hypprior}), with $p_1 \equiv 0.5$ (i.e.\ we have no prior opinion as to whether the driver has been engaging in undeclared charging during the certified interval). If $H_1$ is realized, then we draw from \eqref{eq:xUH1}, with $\overline{x}_U \equiv E_{max}$ (Table~\ref{tab:vehicle_params_specifc}). If $H_0$, then we set $x_U=0$. We undertake 10,000 such trials for each of $y\in\{s,w\}$ (Section~\ref{sec:auxpower}), so that approximately 5,000 are $H_1$ trials, and 5,000 are $H_0$ trials, in each season. In all these trials, we hold  the traffic (${\mathbf d}_G$) and  EV ($\theta$) operating conditions constant, and so the seasonally-adjusted empirical distributions, $\tilde{\mathsf{F}}(x_C \mid  {\mathbf d}_G, \theta, y)$,  of cumulative energy drawn from the battery, $x_C$, during the (multi-trip) certified interval is invariant across the trials, being those adopted in the scenarios above (Figure~\ref{fig:PerformanceNoCheating}).    

\begin{figure}[H]
    \centering
    \subfigure[]{\includegraphics[width=0.45\linewidth]{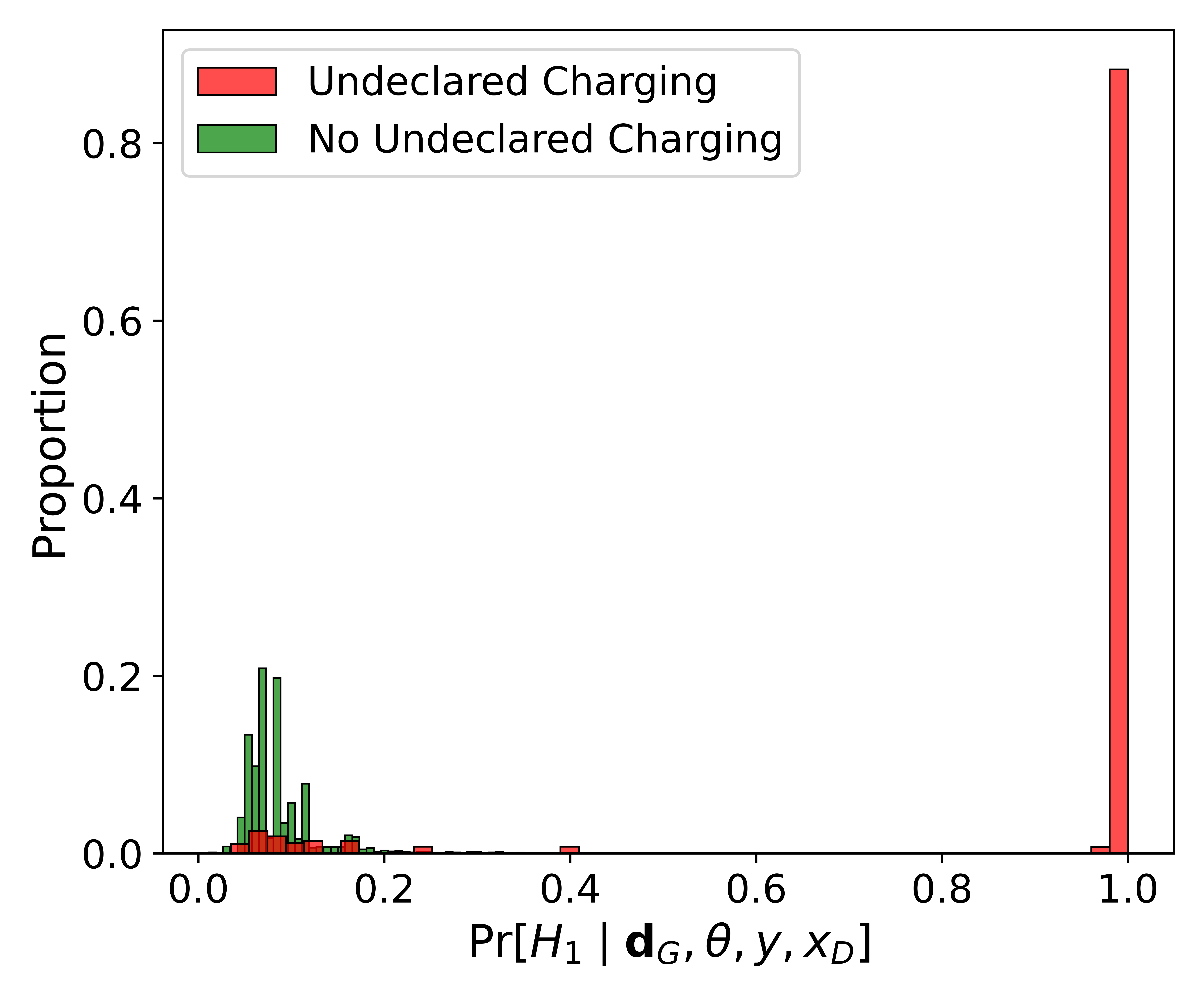}}
    \subfigure[]{\includegraphics[width=0.45\linewidth]{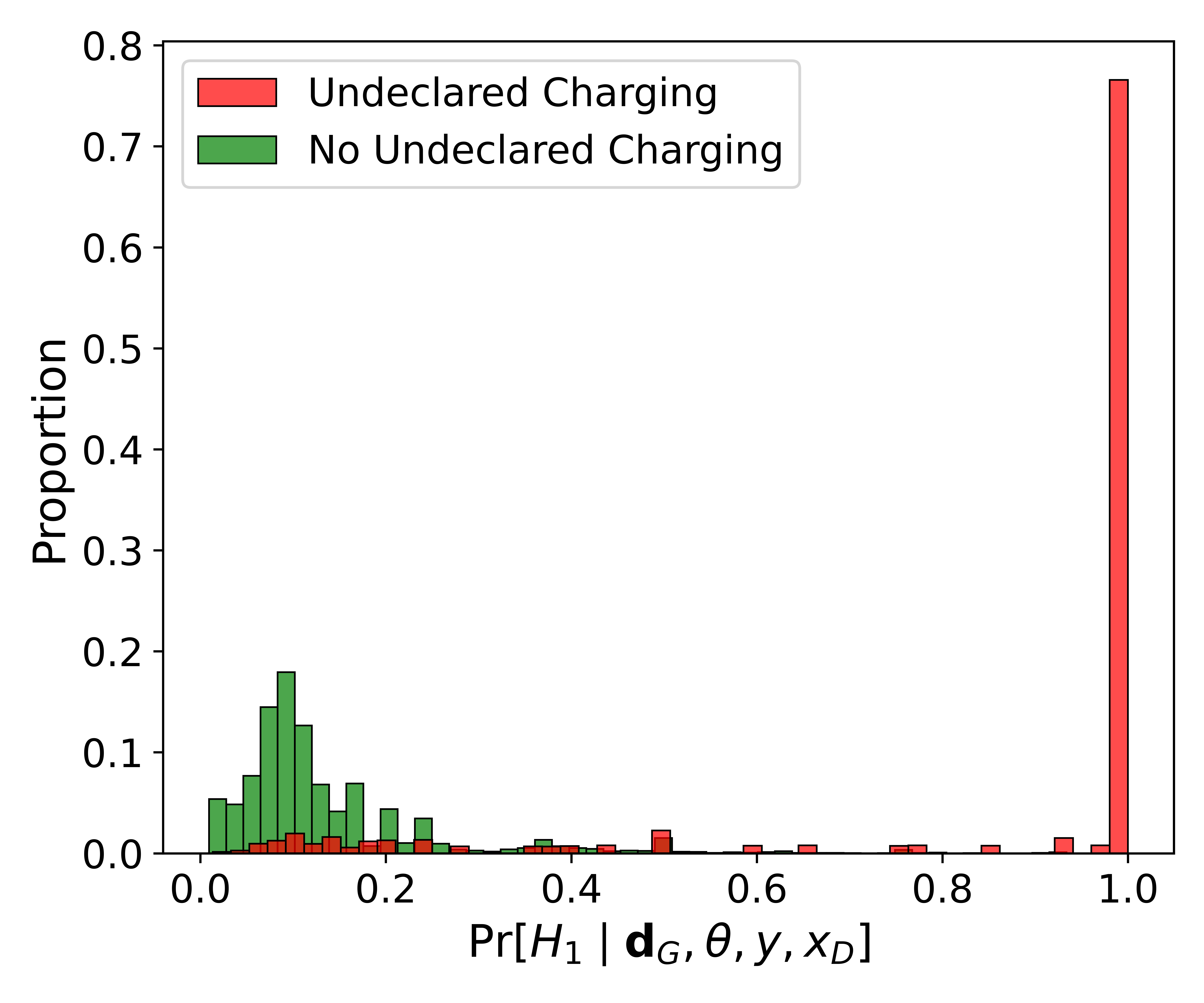}}
    \caption{Histograms of 10,000 trials, in (a) summer, and (b) winter, of the Bayesian hypothesis test for undeclared charging, illustrating $\Pr[H_1 \mid {\mathbf d}_G, \theta, y, x_D]$ in each case. These are segmented into the cases in which  there has, indeed, been undeclared charging during the certified interval (i.e.\ $H_1$, in red), and those in which there was no undeclared  charging  ($H_0$, in green). }
    \label{fig:Cheating_Identification}
\end{figure}

Recall that this is a Bayesian hypothesis test, yielding the {\em a posteriori\/} probability,  $\Pr[H_1 \mid {\mathbf d}_G, \theta, y, x_D]$, as its natural output. We assess the performance of the test via the empirical distributions (histograms) of  $\Pr[H_1 \mid \cdot, x_D]$ for the two prior cases of the driver ($H_1$ and $H_0$), and for the two seasons (Figure~\ref{fig:Cheating_Identification}). 

Clearly, the test performs excellently in discriminating between the two cases, $H_1$ (undeclared charging took place during the certified interval) and $H_0$ (there was no undeclared charging). It should be noted that this MC simulation study was performed under conservative---and therefore testing---conditions, in which there is a high probability of behaviour aimed at fooling the algorithm; i.e.\ $H_1=1$ and $x_U \ll E_{max}$ (as studied in scenario 2). For example, $\Pr[x_U < 0.5 E_{max} \mid H_1] \equiv 0.5$ (\ref{eq:xUH1}), owing to the conservative uniform prior elicitation of $\mathsf{F}(x_U \mid H_1)$ (\ref{eq:xUH1}) (see the red histograms in Figure~\ref{fig:Convolution}). 

In the prior-$H_1$ cases (Figure~\ref{fig:Cheating_Identification}), $\Pr[H_1 \mid {\mathbf d}_G, \theta, y, x_D]=1$ with high probability. Nevertheless, we note that---with probabilities of about 0.2---the test computes $\Pr[H_1 \mid \cdot, x_D]\ll 1$. This is a consequence of the case demonstrated in scenario~2, in which a driver typically attempts to fool the test by charging to a low proportion of $E_{max}$ during the certified interval. 

In the prior $H_0$ cases, there is no confounding with $H_1$ (i.e.\ $\Pr[H_1 \mid {\mathbf d}_G, \theta, y, x_D]\ll 1$ always), but there is more statistical uncertainty in these probabilities. This is a function of the intrinsic (and modelled) uncertainty in the passenger-loading of the EV (\ref{eq:GMM}) and in the driver's use of the auxiliary devices (\ref{eq:Gam}), particularly in winter (Figure~\ref{fig:Gam}).

\begin{figure}[H]
    \centering
    \subfigure[]{\includegraphics[width=0.45\linewidth]{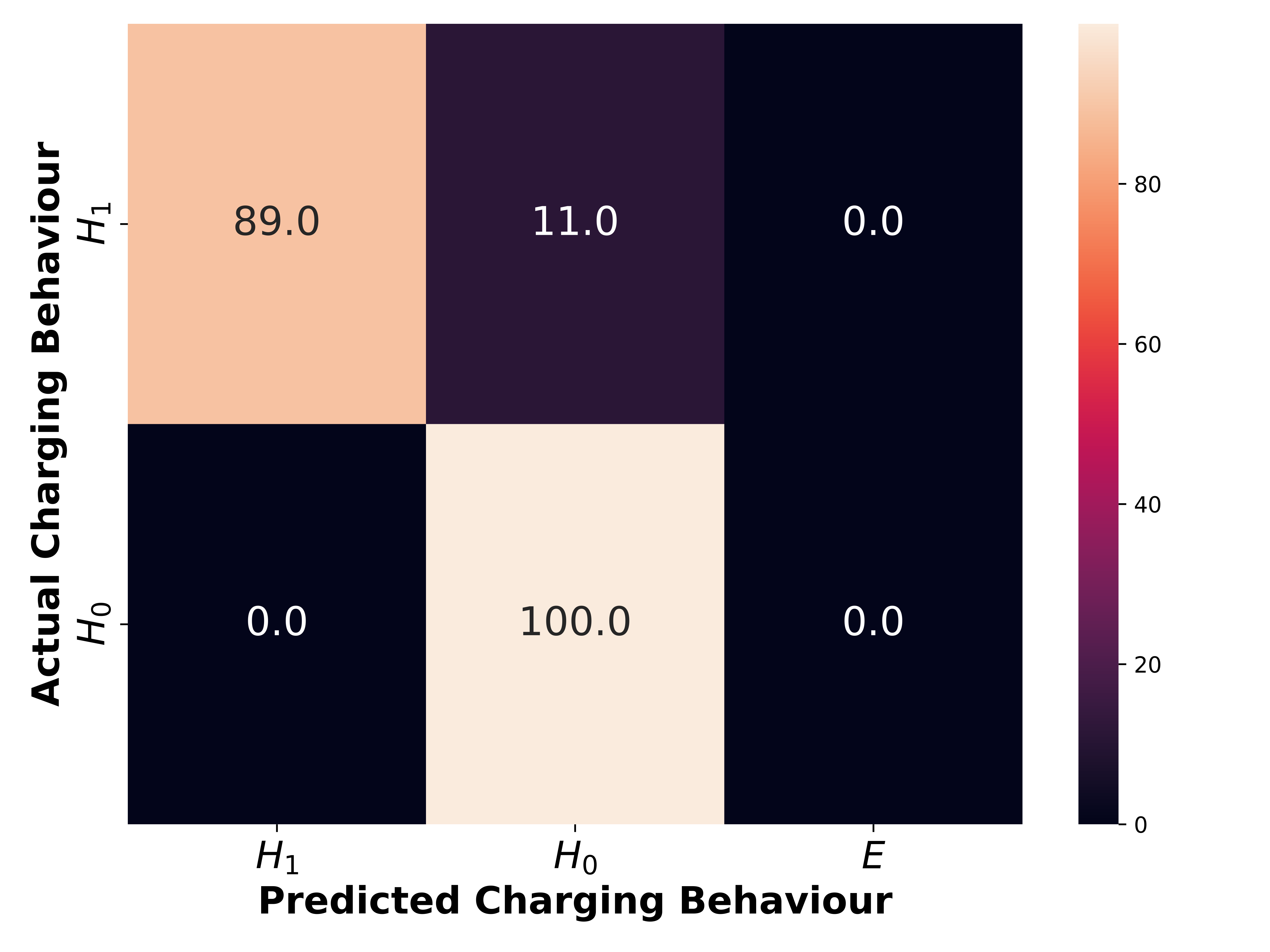}}
    \subfigure[]{\includegraphics[width=0.45\linewidth]{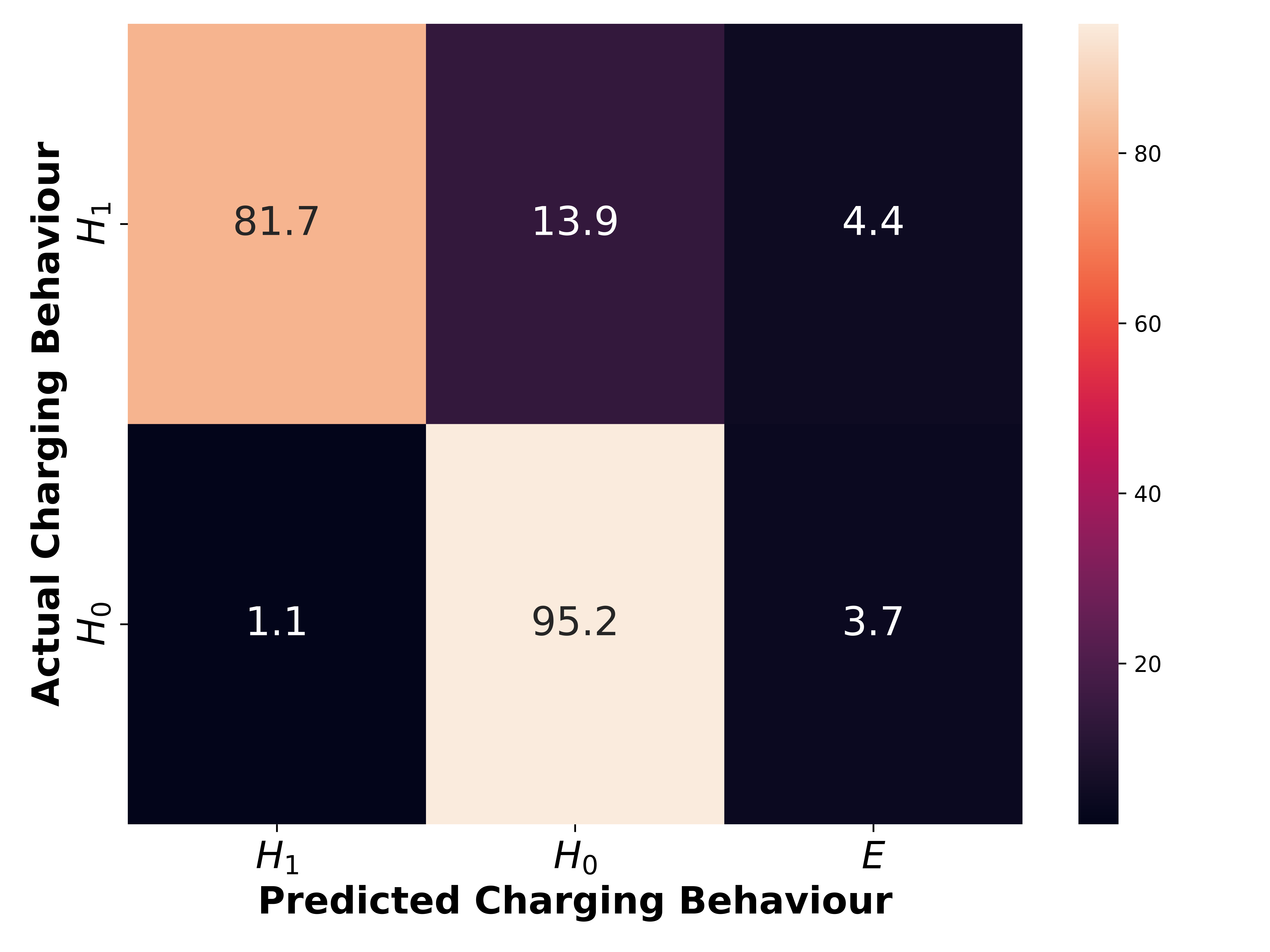}}
    \caption{Confusion matrices of the Bayesian hypothesis test in (a) summer, and (b) winter, when $x_U\mid H_1 \in (0,E_{max})$, as in Figure~\ref{fig:Cheating_Identification}. The matrices record the ternary decision (\ref{eq:eras}) percentages for both types of driver behaviour.}
    \label{fig:Confusion_Matrix}
\end{figure}

If a hard $H_1$-vs-$H_0$ decision is to be implemented, then this is done by thresholding the {\em a posteriori\/} probabilities, $\Pr[H_1 \mid {\mathbf d}_G, \theta, y, x_D]$ \eqref{eq:Bayes}, in an appropriate manner. We censor equivocal inferences by designing a ternary decision, involving an erasure (i.e.\ censored) state, $E$, in which case no decision about the behaviour of the driver is made: 
\begin{equation}
 \Pr[H_1 \mid {\mathbf d}_G, \theta, y, x_D]  \in\left\{ 
\begin{array}{rll}
[0,0.4]&\Rightarrow &H_0,\\
&&\\
(0.6,1.0]&\Rightarrow &H_1,\\
&&\\
\mathsf{otherwise} &\Rightarrow &E. 
\end{array}
\right.
\label{eq:eras}
\end{equation}
We implement  this for each of the 10,000 trials (for each season) of the MC simulation (Figure~\ref{fig:Cheating_Identification}), yielding the seasonally dependent confusion matrices in Figure~\ref{fig:Confusion_Matrix}. Conditional on erasure of  the $E$ states, 
these statistics yield test sensitivities of $89\%$ in summer and $85.8\%$ in winter. The specificities of the test are $100\%$ and $98.7\%$, respectively.

\begin{figure}[H]
    \centering
    \subfigure[]{\includegraphics[width=0.45\linewidth]{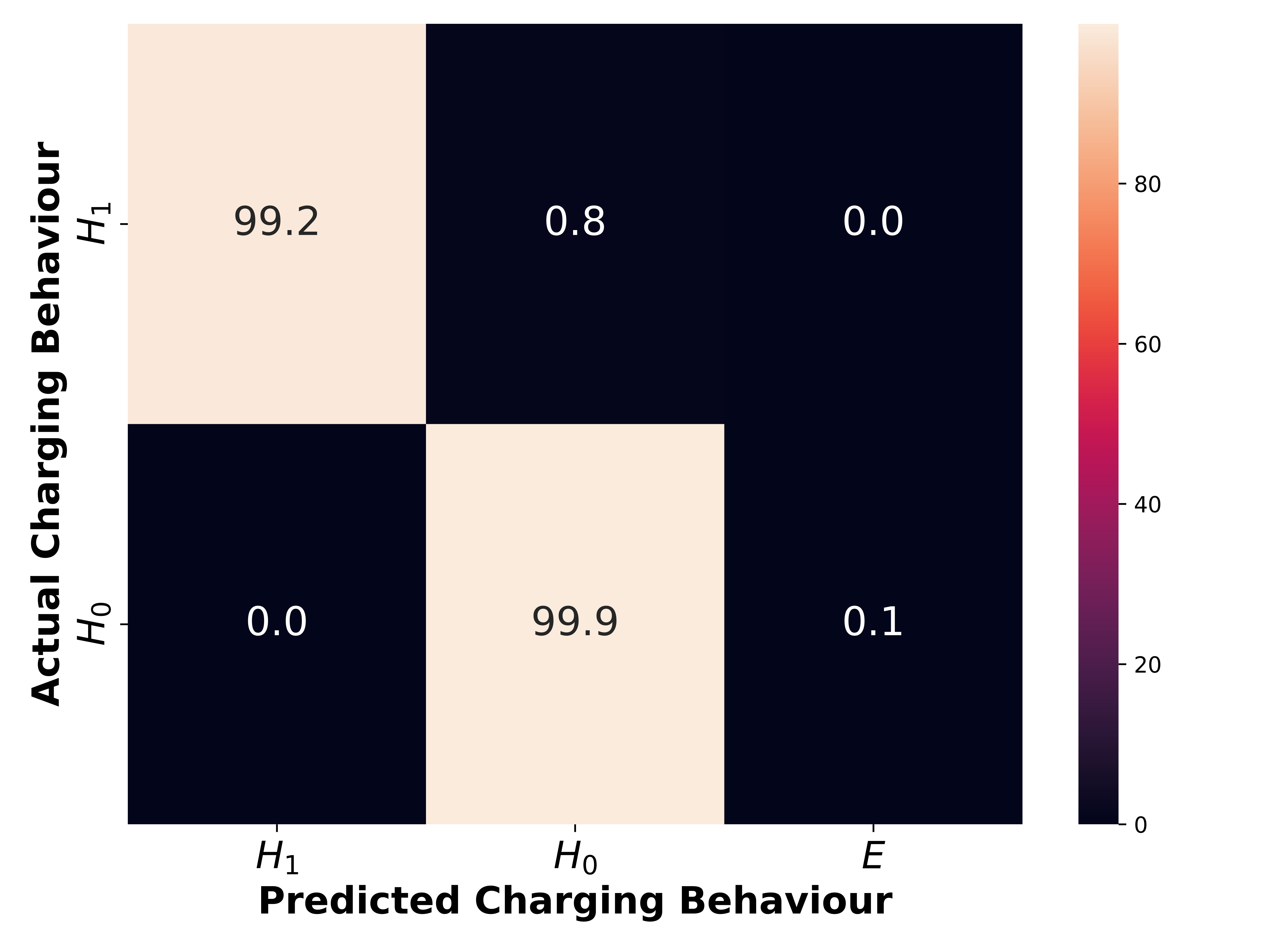}}
    \subfigure[]{\includegraphics[width=0.45\linewidth]{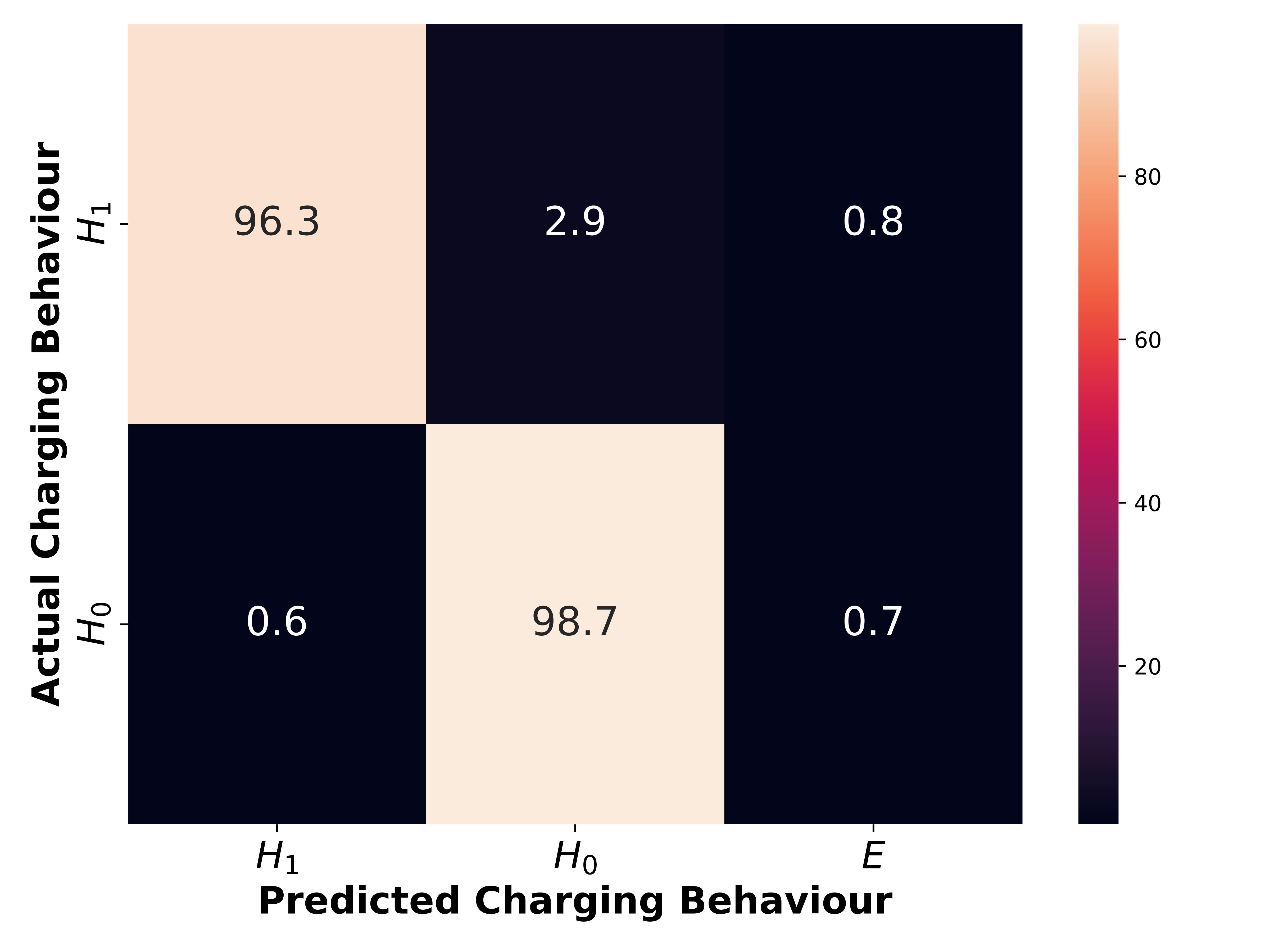}}
    \caption{Confusion matrices of the Bayesian hypothesis test in (a) summer, and (b) winter, when $x_U\mid H_1 \in (0.2 E_{max},E_{max})$. The matrices again record the ternary decision (\ref{eq:eras}) percentages for both types of driver behaviour.}
    \label{fig:Confusion_Matrix20-100}
\end{figure}

In order to validate the proposition, above, regarding the conservative (i.e.\ stressful) operating conditions which have been implemented in our MC simulations (Figures~\ref{fig:Cheating_Identification} and \ref{fig:Confusion_Matrix})---where 50\% of cases of undetected charging are `small' (i.e. $x_U \leq 0.5E_{max}$)---we adjust the range of $x_U\mid H_1$ in \eqref{eq:xUH1} to $(0.2E_{max},E_{max}]$, reducing the probability that undeclared charging is small to $38\%$. As predicted, the performance of our Bayesian hypothesis test of undeclared charging is then even better (Figure~\ref{fig:Confusion_Matrix20-100}), 
with the sensitivities rising to $99.2\%$ and $97\%$, respectively, while the specificities rise to $100\%$ and $99.4\%$, respectively.

\section{Discussion and concluding remarks}
\label{sec:disc}
There are several benefits attached to the  Bayesian formulation \eqref{eq:Bayes} of the test for undeclared charging:
\begin{itemize}
    \item[(i)] 
    The algorithm calculates the EV's  {\em a posteriori\/} probability \eqref{eq:Bayes}, and provides it to the operator of the green-energy certification scheme,
    as soon as it is connected to the next charging station belonging to the scheme. As we have seen \eqref{eq:eras}, this can be used to make a decision as to whether the EV was subject to undeclared---and therefore possibly not green-energy---charging since the last certified charging event. This decision can be used to apply or withhold either an incentive/bonus or a penalty/levy applied to the driver by the operator of the scheme. However, a fairer approach is to use this probability to {\em weight\/} the bonus or levy. For instance, in an incentivization scheme, a simple approach would be to apply the bonus, $(1-\Pr[H_1 \mid \cdot, x_D])G$, where $G$ is the maximal bonus for compliance with the scheme---when the EV is charged at a certified charging station.
    \item[(ii)] The sequential nature of Bayesian inference allows probabilities to be propagated from one certified interval to the next. The {\em a posteriori\/} probability, $\Pr[H_1 \mid \cdot, x_D]$ (\ref{eq:Bayes}), can be assigned as the {\em priori\/} probability, $p_1 \equiv \Pr[H_1]$ \eqref{eq:hypprior}, for the next charging interval. Again, a fairer scheme would seek to model non-stationarity in the driver's behaviour over time. One approach involves the estimation of a forgetting factor, $0\leq \lambda \leq 1$, based on an assessment of the pseudo-stationary window-length of the driver's behaviour (high $\lambda$ corresponding to slow change in behaviour)~\cite{forget}, so that
    \[
    p_{1, i+1} \equiv \lambda \Pr  \hspace*{-.15cm}\;_i [H_1 \mid \cdot, x_D],
    \]
    where $i$ indexes the sequence of certified charging intervals. There is more work to be done in modelling nonstationary dynamics in the driver's behaviour with respect to the certification scheme. Indeed, the same sequential inference procedure could also be used to sharpen the undeclared charging prior,  $\mathsf{F}(x_U \mid H_1)$ \eqref{eq:xUH1}, by updating it to the {\em a posteriori\/} distribution, $\mathsf{F}(x_U \mid {\mathbf d}_G, \theta, y, x_D, H_1)$, based on the evidence from the last certified interval. This would then be propagated forward (perhaps with forgetting) as the $x_U \mid H_1$ prior \eqref{eq:xUH1} for the next interval.  
    \item[(iii)] The parametric Bayesian framework of this paper requires all stochastic modelling assumptions,  and the knowledge base, $\{{\mathbf d}_G, \theta, y, x_D\}$, to be explicitly declared. This facilitates the potential inclusion of other factors into the methodology, such as the effect of per-unit green-energy price volatility on $x_U$, which can be processed into the prior, $\mathsf{F}(x_U \mid H_1)$ (\ref{eq:xUH1}). In this paper, we have adopted noninformative or conservative prior parameters in several instances, such as in $\mathsf{F}(m_{peop})$ (\ref{eq:GMM})  and $\mathsf{F}(W\mid y)$ (\ref{eq:Gam}), as well as in $\mathsf{F}(x_U \mid H_1)$ (\ref{eq:xUH1}). These could all be further data-informed by recourse to publicly available---or driver-supplied---knowledge.
\end{itemize}

On this last point, we have noted (in scenario~2 in  Section~\ref{sec:scen3}) that low-$x_U$ undeclared charging can lead to false negatives in a small number of cases (see Figure~\ref{fig:Cheating_Identification}(b), and the top row of the confusion  matrix in Figure~\ref{fig:Confusion_Matrix}(b)). We  see different---potentially deleterious---effects of uncertainty  in the heavy right-tails of the $\Pr[H_1 \mid \cdot, x_D]$ distributions  in the $H_0$ cases (Figure~\ref{fig:Cheating_Identification}). While these tails do not incur a significant rate of false-positives, the problem is again  caused by the significant uncertainty engendered in the consumed energy predictor, $\tilde{\mathsf{F}}(x_C \mid  {\mathbf d}_G, \theta, y)$ \eqref{eq:xCdist}, particularly because of the uncertainties in auxiliary power, $W$ \eqref{eq:Gam}. In common with {\em loyalty schemes\/}, widely adopted in retail, our framework can be adapted to provide a further bonus to drivers who cooperate in sharing data that improve the statistical performance of the predictor. As ever, the cost to the driver---in terms of loss of privacy---needs to assessed in any fair implementation of the scheme that relies on private-data-driven elicitation of the priors.



Finally, in order to support the incentivization scheme based on our detection algorithm, it would be necessary to implement it on a {\em distributed ledger technology (DLT)\/} platform,  such as blockchain. This would  enable secure, transparent, and tamper-proof recording of charging events by drivers and their energy consumption data. Smart contracts could be developed, which would  apply rewards automatically, based on processing verified charging data,  thereby enhancing trust and reducing administrative effort.

\subsection{Concluding remarks}
\label{sec:conc}

The  paper has developed a novel methodology for Bayesian detection of undeclared charging events for use in a green-energy certification scheme. Simulated GPS data, EV parameters and the measured differential SoC have been processed into a probability of undeclared charging during intervals between certified charging events. Extensive simulations have demonstrated the excellent  statistical performance of the approach under realistic operating conditions. We have outlined how the Bayesian formulation of the test can serve the needs of an incentivization scheme for green-energy consumption, and how it can  readily exploit other knowledge sources.

\bibliographystyle{plain}
\bibliography{PaperBibliography}

\end{document}